\documentclass{amsart}
\usepackage{amscd,amssymb}
\usepackage[mathcal]{euscript}
\newtheorem{theorem}{Theorem}[subsection]
\newtheorem{lemma}[theorem]{Lemma}
\newtheorem{proposition}[theorem]{Proposition}

\newtheorem{corollary}[theorem]{Corollary}
\theoremstyle{definition}
\newtheorem{definition}[theorem]{Definition}
\newtheorem{remark}[theorem]{Remark}

\numberwithin{equation}{subsection}

\DeclareMathOperator{\Hom}{Hom}

\DeclareMathOperator{\map}{map}

\DeclareMathOperator{\sh}{sh}
\DeclareMathOperator{\thh}{THH}
\DeclareMathOperator{\THH}{THH}
\DeclareMathOperator{\tk}{thh^k}

\DeclareMathOperator{\h}{T}
\DeclareMathOperator{\tkt}{tHH^k}

\DeclareMathOperator{\tHH}{tHH}
\DeclareMathOperator{\hocolim}{hocolim}
\DeclareMathOperator{\colim}{colim}
\DeclareMathOperator{\Salg}{S-alg}
\newcommand{\parallelarrows}[1]{\begin{array}{c} {\hbox to
#1{\rightarrowfill}}  \vspace{-0.35cm} \\ {\hbox to
#1{\rightarrowfill}} \end{array}}

\newcommand{\spec}{Sp^{\Sigma }}
\newcommand{\BF}{Sp^{\mathbf{N}}}
\newcommand{\point}{{\textstyle {\bf \cdot}}}
\newcommand{\iso}{\cong}
\newcommand{\hoeq}{\simeq}
\newcommand{\K}{[r]_+}
\newcommand{\p}{$\pi_*$}
\newcommand{\spectra}{\spec}

\newcommand{\sset}{\mathcal{S}_{*}}
\newcommand{\ssets}{\sset}
\newcommand{\tensor}{\otimes}
\newcommand{\otime}{\wedge}
\newcommand{\sm}{\wedge}
\newcommand{\sms}{\wedge_S}

\newcommand{\smasRe}{\wedge_{R^e}}
\newcommand{\smasR}{\wedge_{R}}
\newcommand{\smashR}{\wedge_{R}}
\newcommand{\smashs}{\wedge_S}
\newcommand{\Smash}{\wedge_S}
\newcommand{\smashk}{\wedge_k}
\newcommand{\s}{\colon\,}
\newcommand{\N}{\mathbb{N}}
\newcommand{\D}{\mathcal{D}}
\newcommand{\F}{\mathcal{F}}
\newcommand{\id}{{\mbox id}}
\newcommand{\tw}{{\mbox tw}}

\begin{document}

\title[Ring spectra and THH]{Symmetric ring spectra and topological Hochschild homology} 

\date{\today}
\author{Brooke Shipley}
\thanks{Research partially supported by an NSF Postdoctoral Fellowship} 
\address{Department of Mathematics \\ University of Chicago\\ Chicago, IL 60637
\\ USA} 
\email{bshipley@math.uchicago.edu}
\maketitle

\section{Introduction}\label{intro}

The category of symmetric spectra introduced by Jeff Smith is a 
closed symmetric monoidal category whose associated homotopy category
is equivalent to the traditional stable homotopy category, see [HSS].
In this paper, we study symmetric ring spectra, 
{\it i.e.,} the monoids in the category of symmetric spectra. 
The category of
symmetric ring spectra is closely related to the category of ``functors
with smash product defined on spheres" defined for instance in [HM, 2.7].
Actually, the category of symmetric ring spectra is equivalent to the
category of FSPs defined on spheres if the usual connectivity and
convergence conditions on FSPs are removed. 

The choices of equivalences must also be changed when considering
symmetric ring spectra instead of FSPs on spheres.  A map of FSPs on spheres
is a weak equivalence when the map is a \p-isomorphism, {\it i.e.,} when 
it induces an isomorphism in the associated stable homotopy groups.    
As with symmetric spectra, one must consider a broader class
of equivalences called stable equivalences when working with
symmetric ring spectra, see \ref{stable eq.}.
In section \ref{model}, the model category
structure on symmetric ring spectra is defined with these
stable equivalences.  In [MMSS], we show that the associated homotopy 
category is equivalent to the traditional category of $A_\infty$-ring spectra.

Because there are more stable equivalences than \p-isomorphisms, the
classically defined stable homotopy groups are not invariants of the homotopy 
types of symmetric spectra or symmetric ring spectra.  Hence stable
equivalences can be hard to identify.  To remedy this we consider
a detection functor, $D$, which turns stable equivalences into 
\p-isomorphisms.  Theorem \ref{detection-thm}, shows that  
$X \xrightarrow{} Y$ is a stable equivalence if and only if
$DX \xrightarrow{} DY$ is a \p-isomorphism.  Thus, the classical stable
homotopy groups of $DX$ are invariants of the stable homotopy type of $X$.
There is also a spectral sequence for calculating the classical stable homotopy
groups of $DX$, see Proposition \ref{spec.seq}.

The category of FSPs was defined in [B] in order to define the topological
Hochschild homology for an associative ring spectrum $R$.  In section
\ref{thh.1}, three different definitions of
topological Hochschild homology for a symmetric ring spectrum are considered;
B\"okstedt's original definition restated for symmetric ring spectra, 
\ref{thhr}, a derived smash product definition, \ref{tk}, and a 
definition which mimics
the standard Hochschild complex from algebra, \ref{tkt}.    
Theorems~\ref{thh comparison} and \ref{thm-thh-many} 
show that under certain cofibrancy conditions these definitions
all agree.  

Perhaps the most surprising of these results is the agreement of 
B\"okstedt's definition with the others without any connectivity or convergence
conditions.  Some conditions are indeed 
necessary to apply B\"okstedt's approximation theorem, [B, 1.6], though 
the usual connectivity and convergence conditions can be weakened, see 
Corollary~\ref{cor-sf}. 
For spectra which do not satisfy these hypotheses
the model category structure on symmetric ring spectra is used instead to prove 
comparison results such as Theorem \ref{detection-thm}.  
Also, without
any extra conditions, B\"okstedt's original definition of THH takes stable
equivalences of symmetric ring spectra to \p-isomorphisms, see 
Corollary \ref{thh-stable}.  See also Remark \ref{FSP}.

{\em Outline.}
In the first section we recall various definitions from [HSS], define
symmetric ring spectra, and discuss homotopy colimits.  
Section \ref{model} uses the results of [SS] and [HSS] to establish 
model category structures for symmetric ring spectra,
for $R$-modules over any symmetric ring spectrum $R$, and for $R$-algebras over
any commutative symmetric ring spectrum $R$.  
In section \ref{hocolim} we define the 
homotopy colimit of diagrams of symmetric spectra and state several comparison 
results for homotopy colimits which are used in sections \ref{D} and
\ref{thh.1}.   A functor, $D$, which detects stable equivalences is defined in
section \ref{D}.  In section \ref{thh.1} three 
different definitions of topological Hochschild homology are defined and 
compared.    
The detection functor from section \ref{D} is used in one of the comparisons in 
section \ref{thh.1}.  

{\em Acknowledgments.}
I would like to thank Bill Dwyer and Jeff Smith for many helpful conversations 
throughout this work.  I also benefited from discussions with Lars Hesselholt, Mark Hovey, Mike Mandell, Haynes Miller, Charles Rezk, and 
Stefan Schwede.    

\section{Basic definitions}\label{basics}
In this section we state the basic definitions which are needed in
sections \ref{D} and \ref{thh.1}. 
Most of the definitions in the first subsection come from 
[HSS].
The second subsection, \ref{sec-srs}, considers symmetric ring spectra and
model categories for $R$-modules and $R$-algebras.  In the last
subsection, \ref{hocolim}, we consider the properties of the homotopy
colimit needed for sections \ref{D} and \ref{thh.1}.

\subsection{Symmetric spectra}\label{sec-spec}
We first define the symmetric monoidal category of symmetric spectra.  
Next we define certain model category structures on symmetric spectra.  
Then we consider a subcategory of symmetric spectra, the semistable 
spectra, between which stable equivalences are exactly the \p-isomorphisms.
Throughout this paper ``space" means simplicial set, except in Remark 
\ref{top-spaces}.

\begin{definition}\label{def-symmetric-sequences}
Let $\Sigma$ be the skeleton of the category of finite sets and bijections
with objects ${\bf n}=\{1,\cdots, n\}$.
The category of {\em symmetric sequences} $\sset^{\Sigma}$ is the
category of functors from $\Sigma $ to $\ssets$, the category of pointed
simplicial sets.  
Thus, a symmetric sequence is a sequence $X_{n}$ of spaces, where
$X_{n}$ is equipped with a basepoint preserving
action of the symmetric group $\Sigma _{n}$. 
\end{definition}

This category of symmetric sequences is a symmetric monoidal category with
the following definition of the tensor product of two symmetric sequences.

\begin{definition}[HSS]\label{def-tensor-sequences}
Given two symmetric sequences $X$ and $Y$, we define their tensor
product, $X\otimes Y$, 
\[
(X\otimes Y)_{n}=\bigvee _{p+q=n} \Sigma _{n}^{+}\wedge _{\Sigma
_{p}\times \Sigma _{q}} (X_{p}\wedge Y_{q}).
\]
Note here, as elsewhere in the paper, we sometimes denote an extra
basepoint as $X^+$ to make the notation more readable.  
\end{definition}

Let $S^1= \Delta[1]/\dot{\Delta}[1]$ and $S^n = (S^1)^{\wedge n}$ for $n > 1$.  
Then $S=(S^0, S^1, \cdots, S^n, \cdots)$ is a symmetric sequence.  
In fact, $S$ is a commutative monoid.  

\begin{definition} The category of {\em symmetric spectra}, $\spec$, is the
category of left $S$-modules in the category of symmetric sequences.
\end{definition}

A symmetric spectrum is then a sequence of
pointed spaces with a left, pointed $\Sigma_p$ action on
$X_p$ and associative, unital, $\Sigma _{p} 
\times \Sigma_{q}$-equivariant maps $S^{p} \wedge X_q \xrightarrow{} X_{p+q}$.

The category of symmetric spectra is a symmetric monoidal category with
the following definition of the smash product of two symmetric spectra,
see [HSS]. 

\begin{definition}
Given two symmetric spectra $X$ and $Y$ we define their smash product
$X \wedge_S Y$ as the coequalizer of the two maps 
$$X \tensor S \tensor Y \parallelarrows{0.5cm} 
X \tensor Y.$$
\end{definition}

We now describe certain symmetric spectra which play an important role in the 
model category structures and in the later sections of this paper.  Let
$I$ be the skeleton of the category of finite sets and injections  
with objects ${\bf n}$.  Note that
$\hom_I(n,m) \iso \Sigma_m / \Sigma_{m-n}$ as $\Sigma_m$ sets. 

\begin{definition}\label{free}
Define $F_n\s \sset \to \spec$ by $F_nK=S \otimes G_nK$ where
$G_nK$ is the symmetric sequence with $\hom_I(n,n)_+ \sm K$ in degree $n$
and the basepoint elsewhere.  So $(F_nK)_m= \Sigma_m^+ \sm_{\Sigma_{m-n}} 
S^{m-n}\sm K\hoeq\hom_I(n,m)_+ \sm S^{m-n} \sm K$ where $S^n=*$ for $n < 0$.  
\end{definition} 

$F_n$ is left adjoint to the $n$th evaluation functor $Ev_n: \spec
\xrightarrow{} \ssets$ where $Ev_n(X)=X_n$.  
There is a natural isomorphism $F_nK \smashs F_mL\xrightarrow{}
F_{n+m}(K \wedge L)$.

\emph{Model category structures.}
There are two model category structures on symmetric spectra  
which we consider; the injective model category and the stable model
category.  The injective model category is a stepping stone for defining
the stable model category.  The homotopy category associated to
the stable model category is equivalent to the stable homotopy category of
spectra, see [HSS]. Hence, the stable model category is the model category
which we refer to most often. 
See [Q] or [DS] for the basic definitions for model categories.

\begin{definition} Let $f:X \xrightarrow{} Y$ be a map in $\spec$.
The map $f$ is a {\em level equivalence} if each $f_n:X_n \xrightarrow{}
Y_n$ is a weak equivalence of spaces, ignoring the $\Sigma_n$ action.
It is a {\em level cofibration} if each $f_n$ is a cofibration of
spaces.  
\end{definition}

With level equivalences,
level cofibrations, and fibrations the maps with the
right lifting property with respect to all maps which are trivial
cofibrations, $\spec$ forms a simplicial model category
referred to as the {\em injective model category}.  A cofibration
here is just a monomorphism.
A fibrant object here is called an \emph{injective spectrum}.  An injective
spectrum is a  
spectrum with the extension property with respect to every monomorphism that
is a level equivalence. 

To define the equivalences for the stable model category we first need the
following definition.

\begin{definition} A spectrum $X$ in $\spec$ is an {\em $\Omega$-spectrum}
if $X$ is fibrant on each level and the adjoint to the 
structure map $S^1 \wedge X_n \xrightarrow{} X_{n+1}$
is a weak equivalence of spaces for each $n$.
\end{definition}

Define shifting down functors, $\sh_n:\spec \xrightarrow{} \spec$, by 
$\sh_n(X)_k= X_{n+k}$. Then the adjoint of the structure maps
give a map $i_X \s X \to \Omega \sh_1 X$.    

\begin{lemma} $X$ is an $\Omega$-spectrum if and only if 
$X$ is level fibrant and $X \to \Omega sh_1 X$ is a level equivalence. 
\end{lemma}

\begin{definition}\label{stable eq.}
Let $f: X \xrightarrow{} Y$ be a map in $\spec$.
The map $f:X \xrightarrow{} Y$ is a {\em stable equivalence}
if 
$$\pi_0\map(Y,Z) \xrightarrow{} \pi_0\map(X,Z)$$ 
is an isomorphism for all injective $\Omega$-spectra $Z$.  The map $f$ is a
{\em stable cofibration} if it has the left lifting property with respect
to each level trivial fibration, {\it i.e.,} a map that is a trivial fibration on 
each level.  The map $f$
is a {\em stable fibration} if it has the right lifting property with
respect to each map which is both a stable cofibration and a 
stable equivalence.
\end{definition}

\begin{theorem}[HSS]\label{stable model cat.}
With these definitions of stable equivalences, stable cofibrations, and
stable fibrations, $\spec$ forms a model category referred to as
the {\em stable model category}.  A  map is a stable trivial fibration 
if and only if it is a level trivial fibration.  Moreover, the
fibrant objects are the $\Omega$-spectra and a map between $\Omega$-spectra
is a stable equivalence if and only if it is a level equivalence.
\end{theorem} 

As shown in [HSS], the stable model category is in fact a cofibrantly 
generated model
category.  In particular, this means that a transfinite version of
Quillen's small object argument, [Q, II 3.4], exists.  This argument
is central to the proofs in sections \ref{D} and \ref{thh.1}. 

\begin{proposition}[HSS]\label{cof-gen}
There is a set of maps $J$ in $\spec$ such that any stable trivial
cofibration is a retract of a directed colimit of pushouts
of maps in $J$.  Similarly, there is a set of maps $I$ which generate
the stable cofibrations.   These maps are called the generating stable
(trivial) cofibrations.
\end{proposition} 

Let $I=\{F_n(\dot{\Delta}[k]_+) \xrightarrow{}
F_n({\Delta}[k]_+)\}$ and $J'= \{F_n(\Delta^l [k]_+) \xrightarrow{}
F_n(\Delta [k]_+)\}$.  $I$ is the set mentioned in the proposition which
generates the stable cofibrations.  The set $J$ in the proposition is
the union of $J'$ with a set $K$ which we now describe.   

First, consider the map $\sigma: F_1(S^1) \xrightarrow{}
F_0(S^0)$ which is adjoint to the identity map $S^1 \xrightarrow{} S^1=
Ev_1(F_0S^0)$.   There is a factorization of $\sigma$ as a level cofibration
$c:F_1(S^1) \xrightarrow{} C$ followed by a level trivial fibration 
$r:C \xrightarrow{} F_0(S^0)$.  
$C$ is defined as the pushout in the following square.
\[
\begin{CD}
F_1S^1 \sms F_0(\Delta[0]_+) @>{\sigma}>> F_0S^0\\
@V{1\sms {i_0}}VV   @ViVV\\
F_1S^1 \sms F_0(\Delta[1]_+) @>b>> C
\end{CD}
\]
Define $c:F_1S^1 \xrightarrow{} C$ as the composite $b \circ (1\sms i_1)$.
In [HSS], $c$ is shown to be a stable trivial cofibration.
The left inverse, $r_0$, to $i_0$ induces a map 
$\sigma \circ (1\sms r_0):F_1S^1\sms F_0(\Delta[1]_+) \xrightarrow{} F_0S^0$.  
Use this map and the identity map on $F_0S^0$ to induce 
$r\s C \to F_0S^0$ using the property of the pushout. 

The generating stable trivial cofibrations are built from the map 
$c\s F_1S^1 \to C$ as follows.  
Let $P_{m,r}$ be the pushout of the following square.
\[
\begin{CD}
F_1S^1 \smashs F_m(\dot{\Delta}[r]_+) @>{c \sms 1}>> C \smashs
F_m(\dot{\Delta}[r]_+)\\
@V{1 \sms F_mg_r}VV   @VVV\\
F_1S^1 \smashs F_m(\Delta[r]_+) @>>> P_{m,r} 
\end{CD}
\]
The map $P(c, F_mg_r): P_{m,r} \xrightarrow{} C \smashs F_m(\Delta[r]_+)$ 
is induced by the property of the pushout by the maps 
$c \sms 1_{F_m(\Delta[r]_+)}$ and $1_C \sms F_mg_r$.
Let $J=J' \cup K$ where $K=\{P(c, F_mg_r), m, r \geq 0\}$  and
$J'$ is the set of generating level trivial cofibrations defined above.  
  
Quillen's small object argument [Q, p.\ II 3.4] has 
an analogue which allows one to functorially factor maps whenever
the model category is cofibrantly generated, see [HSS].   

\begin{definition}\label{L} Let $L$ be the functorial {\em stable fibrant 
replacement}
functor defined by functorially factoring the map $X \xrightarrow{} *$ 
into a stable trivial cofibration, $X \xrightarrow{} LX$ and a stable fibration
$LX \xrightarrow{} *$.  This is the factorization one defines using
the small object argument applied to the set of maps $J$.  Using $J'$ 
instead, one defines $L'$ as the functorial {\em level fibrant replacement}
with $X \to L'X$ a level trivial cofibration and $L'X \to *$ a level fibration.
\end{definition}

This analogue of the small object argument provides a general method
for proving that the class of cofibrations or trivial cofibrations has some 
property.  One only needs to show that the generating maps have some property 
and that the property is preserved under pushouts, colimits, and
retracts.  This general method is used in 
sections \ref{D} and \ref{thh.1}. 

We also need to know that the symmetric monoidal structure and the model
category structure fit together to give a monoidal model category
structure, see [HSS]. 

\begin{proposition}[HSS]\label{monoidal}
The stable model category is a {\em monoidal model category}. 
That is, 
the  symmetric monoidal structure satisfies the following pushout product axiom.
Let $f: A \xrightarrow{} B$ and $g: C \xrightarrow{} D$
be two stable cofibrations.  Then 
$$Q(f,g): (A \wedge_S D) \cup_{A \wedge_S C} (B \wedge_S C) \xrightarrow{} B 
\wedge_S D$$ is a stable cofibration.  
If one of $f$ or $g$ is also a stable equivalence, then so is $Q(f,g)$.   
\end{proposition}

This structure
can also be restricted to a simplicial structure via the functor 
$F_0:\ssets \xrightarrow{} \spec$.  

\emph{Semistable objects and \p-isomorphisms.}\label{subsec-sf}
It is often useful to compare symmetric spectra to the model category
of spectra, $\BF$, defined in [BF].  
There is a forgetful functor $U: \spec \xrightarrow{} \BF$ 
which forgets the action of the symmetric groups and uses the structure
maps $S^1 \wedge X_n \xrightarrow{} X_{1+n}$. 

\begin{definition}\label{shg}
Let $\pi_k(X)= \pi_k(UX)=\colim_i \pi_{k+i} X_i.$
A map $f$ of symmetric spectra is a {\em \p-isomorphism} if it 
induces an isomorphism on these \emph{classical stable homotopy groups}.
\end{definition}

As seen in [HSS] these classical stable homotopy groups are NOT the maps
in the homotopy category of symmetric spectra of the sphere into $X$.
For example $\sigma \s F_1S^1 \xrightarrow{} F_0S^0$ is a stable equivalence
but it is not a \p-isomorphism.  As shown in [HSS], though, a
\p-isomorphism is a particular example of a stable equivalence.
Hence, to avoid confusion, we use the term \p-isomorphism instead
of stable homotopy isomorphism and call these the classical stable homotopy
groups instead of just stable homotopy groups. 
In section~\ref{D} we construct
a functor, $D$, 
which converts any stable equivalence into a \p-isomorphism
between semistable spectra, see Definition~\ref{semistable} below.  

As in [BF], we define a functor $Q$ for symmetric spectra.  

\begin{definition} 
Define $QX= \colim_n\Omega^n L'\sh_n X.$ 
\end{definition}

This functor does not have the same properties as in 
[BF].  For instance $QX$ is not always an $\Omega$-spectrum and $X \xrightarrow{}
QX$ is not always a \p-isomorphism.  One property that does continue to 
hold, however,
is that a map $f$ is a \p-isomorphism if and only if $Qf$ is a level 
equivalence.  Also, $QX$ is always level fibrant. 

\begin{definition}\label{semistable}
A {\em semistable} symmetric spectrum is one
for which the stable fibrant replacement map, $X \xrightarrow{} LX$, is
a \p-isomorphism.
\end{definition}

Of course $X \xrightarrow{} LX$ is always a stable equivalence, but
not all spectra are semistable.  For instance $F_1S^1$ is not
semistable.  Any stably fibrant spectrum, {\it i.e.,} an $\Omega$-spectrum, is  
semistable though.
The following proposition shows that on semistable spectra 
$Q$ has the same properties as in [BF] on $\BF$.     

\begin{proposition}\label{sf}
The following are equivalent.
\begin{enumerate}
\item The symmetric spectrum $X$ is semistable.
\item The map $X \xrightarrow{} \Omega L'\sh_1 X$ is a \p-isomorphism.
\item $X \xrightarrow{} QX$ is a \p-isomorphism.
\item $QX$ is an $\Omega$-spectrum.
\end{enumerate}
\end{proposition}

Before proving this proposition we need the following lemma.

\begin{lemma}\label{lem-spec-Q}
Let $X\in \spec $.  Then 
$\pi _{k}(QX)_{n}$ and $\pi
_{k+1}(QX)_{n+1}$ are isomorphic groups, and 
$i_{QX}\s (QX)_n \to \Omega (Q\sh_1 X)_{n+1}$ 
induces a monomorphism $\pi _{k}(QX)_{n}\xrightarrow{}\pi _{k+1}(QX)_{n+1}$.
\end{lemma}

\begin{proof}
We can assume $X$ is level fibrant.  Both $\pi_{k}(QX)_{n}$ and
$\pi_{k+1}(QX)_{n+1}$ are isomorphic to the $(k-n)$th classical stable
homotopy group of $X$.  However, the map $\pi_{k}i_{QX}$ need not
be an isomorphism.  Indeed, $\pi_{k}i_{QX}$ is the map induced on
the colimit by the vertical maps in the diagram \[
\begin{CD}
\pi _{k}X_{n} @>>> \pi _{k+1}X_{n+1} @>>> \pi _{k+2}X_{n+2} @>>>
\cdots \\ @VVV @VVV @VVV \\ \pi _{k+1}X_{n+1} @>>> \pi
_{k+2}X_{n+2} @>>> \pi _{k+3} X_{n+3} @>>> \cdots
\end{CD}
\]
where the vertical maps are not the same as the horizontal maps,
but differ from them by isomorphisms.  The
induced map on the colimit is injective in such a situation, though
not necessarily surjective.
\end{proof}

\begin{proof}[Proof of Proposition~\ref{sf}]
First we show that (1) implies (2) by using the following diagram.
$$
\begin{CD}
X @>>> \Omega L' \sh_1 X\\
@VVV       @VVV\\
LX @>>> \Omega L' \sh_1 LX
\end{CD}
$$
Since $\Omega L' \sh$ preserves \p-isomorphisms both vertical
arrows are \p-isomorphisms.  The bottom map is a level 
equivalence since $LX$ is an $\Omega$-spectrum.  Hence the top map
is also a \p-isomorphism.

Also, (2) easily implies (3). Since $\Omega$ and $\sh$ commute and both preserve
\p-isomorphisms, $X \to QX$ is a colimit
of \p-isomorphisms provided $X \to \Omega L'\sh_1 X$ is
a \p-isomorphism.

Next we show that (3) and (4) are equivalent.  The map $\pi_*X \to \pi_*QX$ 
factors as $\pi_*X \to \pi_*(QX)_0 \to \pi_*QX$ where the first map here
is an isomorphism by definition.  Then by Lemma~\ref{lem-spec-Q} we see that
$\pi_*(QX)_0 \to$ \p$QX$ is an isomorphism if and only
if $\pi_*(QX)_n \to \pi_{*+1}(QX)_{n+1}$ is an isomorphism for each $n$.  

To see that (3) implies (1), consider the following diagram.
$$
\begin{CD}
X @>>> LX\\
@VVV       @VVV\\
QX @>>> QLX
\end{CD}
$$
By (3) and (4) the left arrow is a \p-isomorphism to an $\Omega$-spectrum, $QX$.
Since $LX$ is an $\Omega$-spectrum the right arrow is a 
level equivalence.  Since $X \to LX$ is a stable equivalence,
the bottom map must also be a stable equivalence.
But a stable equivalence between $\Omega$-spectra
is a level equivalence, so the bottom map is a level equivalence.
Hence the top map is a \p-isomorphism.
\end{proof}

Two classes of semistable spectra are described in the following
proposition.
The second class includes the  connective and convergent spectra. 

\begin{proposition}\label{sf-ex}
\begin{enumerate}
\item If the classical stable homotopy groups of $X$ are all finite then 
$X$ is semistable.
\item Suppose that $X$ is a level fibrant symmetric spectrum 
and there exists some $\alpha > 1$ such that 
$X_{n}\xrightarrow{}\Omega X_{n+1}$ induces an isomorphism 
$\pi_{k}X_{n}\xrightarrow{}\pi _{k+1}X_{n+1}$ for all $k\leq \alpha n$ for
sufficiently large $n$.  Then $X$ is semistable.
\end{enumerate}
\end{proposition}

\begin{proof}
By Lemma~\ref{lem-spec-Q}, $\pi_k(QX)_n \to \pi_{k+1}(QX)_{n+1}$ is a
monomorphism between two groups which are isomorphic. 
In the first case these groups are finite, so this
map must be an isomorphism.  Hence $QX$ is an $\Omega$-spectrum, 
so $X$ is semistable.

For the second part we also show that $QX$ is an $\Omega$-spectrum. 
Since for fixed $k$ the maps $\pi_{k+i}X_{n+i} \to \pi_{k+1+i}X_{n+1+i}$
are isomorphisms for large $i$, $\pi_k(QX)_n \to \pi_{k+1}(QX)_{n+1}$ is
an isomorphism for each $k$ and $n$.
\end{proof}

The next proposition shows that
stable equivalences between semistable spectra are particularly
easy to understand.

\begin{proposition}\label{sf-equivs}
Let $f:X \xrightarrow{} Y$ be a map between two semistable symmetric
spectra.  Then $f$ is a stable equivalence if and only if it is a
\p-isomorphism.
\end{proposition}

\begin{proof}
Since $X \to LX$ and $Y \to LY$ are \p-isomorphisms, $Lf$ is a \p-isomorphism
if an only if $f$ is a \p-isomorphism.  But $Lf$ is a map between 
$\Omega$-spectra so it is a \p-isomorphism if an only if it is a level equivalence.
But, in general, 
$f$ is a stable equivalence if and only if $Lf$ is a level equivalence.
\end{proof}

Finally, we show that any spectrum \p-isomorphic to a semistable spectrum
is itself semistable.

\begin{proposition}\label{sf-p-eq}
If $f\s X \to Y$ is a \p-isomorphism and $Y$ is semistable then $X$ is
semistable.
\end{proposition}

\begin{proof}
Since $Lf$ and $Y \to LY$ are \p-isomorphisms, $X \to LX$ is also a 
\p-isomorphism.
\end{proof}

\subsection{Symmetric ring spectra}\label{sec-srs}
In this section, rings, modules, and algebras are defined for symmetric spectra.
We also discuss the model category structures on these categories. 

\begin{definition}\label{ring} 
A {\em symmetric ring spectrum} is a monoid in the category
of symmetric spectra.  In other words, a symmetric ring spectrum is
a symmetric spectrum, $R$, with maps $\mu: R \smashs R \xrightarrow{} R$ and
$\eta: S \xrightarrow{} R$ such that they are associative and unital,
{\it i.e.,} $\mu \circ ({\mu\Smash \mbox{\id}})= \mu \circ
({\mbox{\id}\Smash \mu}) $ and 
$\mu \circ ({\eta\Smash \mbox{\id}})\iso 
\mbox{\id}\iso
\mu \circ ({\mbox{\id}\Smash \eta})$.
$R$ is called \emph{commutative} if  $\mu \circ {\mbox {\tw}}
= \mu$ where ${\mbox {\tw}}\s R\Smash R \to R\Smash R$ is
the twist isomorphism. 
\end{definition}

Since symmetric ring spectra are the only type of ring spectra
in this paper we also refer to them as simply \emph{ring spectra}.
Using formal properties of symmetric monoidal
categories, one can show that a monoid in the category of $S$-modules 
is the same as a monoid in the category of symmetric sequences
with a monoid map $\eta: S \xrightarrow{} R$ which is central 
in the sense that $\mu \circ ({\mbox{\id}\Smash \eta}) 
\circ {\mbox {\tw}} = \mu \circ ({\eta\Smash \mbox
{\id}})$.

\begin{remark}\label{FSP}
This description of a symmetric ring spectrum 
agrees with the definition of a functor with smash product defined on
spheres as in [HM, 2.7].  The centrality condition mentioned
above, and in [HM, 2.7.ii], is necessary but was not included in
some earlier definitions of FSPs defined on spheres.  
Note, however, that there are no connectivity
(e.g. $F(S^{n+1})$ is $n$-connected) or convergence conditions (e.g. the limit
is attained at a finite stage in the colimit defining $\pi_n$ for each $n$) 
placed on symmetric ring spectra.  These conditions
are usually assumed although not always explicitly stated when using FSPs.  In
particular, these conditions are necessary for applying B\"okstedt's 
approximation theorem, [B, 1.6].  Corollary~\ref{cor-sf} shows that 
a special case
of this approximation theorem holds for any semistable spectrum. 
To consider non-convergent spectra we use 
Theorem~\ref{detection-thm} in place of the 
approximation theorem.  This theorem does not require any 
connectivity or convergence conditions.   
\end{remark}

Proposition~\ref{sf-ex} shows that the connectivity and convergence conditions 
on an FSP ensure that the associated underlying symmetric spectrum is 
semistable.   Proposition~\ref{sf-equivs} shows that stable equivalences 
between such
FSPs are exactly the \p-isomorphisms.  As with the category
of symmetric spectra, inverting the \p-isomorphisms is not enough
to ensure that the homotopy category of symmetric ring
spectra is equivalent to the homotopy category of $A_\infty$-ring spectra.
So once the connectivity and convergence conditions are removed one
must consider stable equivalences instead of just \p-isomorphisms. 

We also need the following definitions of $R$-modules and $R$-algebras 
in later sections.

\begin{definition} Let $R$ be a symmetric ring spectrum.  A (left)
\emph{$R$-module} is a symmetric spectrum $M$ with a map 
$\alpha: R \smashs M \xrightarrow{} M$ that is associative and unital.
\end{definition}

\begin{definition} Let $R$ be a commutative ring spectrum.  
An \emph{$R$-algebra}
is a monoid in the category of $R$-modules.  That is, an $R$-algebra
is a symmetric spectrum $A$ with 
$R$-module maps $\mu:A \smashR A \xrightarrow{} A$ and
$R \xrightarrow{} A$ that satisfy the usual associativity and unity diagrams.
\end{definition} 

Note that symmetric ring spectra are exactly the $S$-algebras. 
 
\emph{Model category structures for rings and modules}\label{model}
We developed techniques in [SS] to form model category structures for 
algebras and modules over a cofibrantly generated, monoidal model category. 
To apply the results from [SS] the smash product is required to satisfy the
monoid axiom, which is verified in [HSS].
Hence, the category of modules over a given symmetric ring spectrum $R$ and 
the category
of $R$-algebras for any given commutative ring spectrum $R$ are model 
categories, see Theorems~\ref{R-mod} and \ref{R-alg}.
Each of these model category structures uses the underlying stable equivalences 
and underlying stable fibrations as the new weak equivalences and fibrations.  
A map is then a cofibration if it has the left lifting property with respect
to the underlying stable trivial fibrations.
These model category structures are used in section~\ref{thh.1}.   

After establishing these model category structures we state certain comparison 
theorems which
show that a weak equivalence of ring spectra induces
an equivalence of the homotopy theories of the respective modules and algebras. 
These comparison theorems show that the cofibrancy condition that appears
in certain theorems in section~\ref{thh.1} is not too restrictive.

In [HSS] $\spec$ is shown to be a cofibrantly generated, monoidal model
category which satisfies the monoid axiom.  Hence Theorem 4.1 in [SS] applies
to give the following results in the category of symmetric spectra.

\begin{theorem}\label{R-mod}
Let $R$ be a symmetric ring spectrum.  Then the category of $R$-modules
is a cofibrantly generated model category with weak equivalences and
fibrations given by the underlying stable model category structure of $\spec$.  
The generating
cofibrations and trivial cofibrations are given by applying $R \smashs -$
to the generating maps in $\spec$.
\end{theorem}

\begin{proof}
This follows from Theorem 4.1 (1) in [SS], once we note that the domains of
the generating cofibrations and trivial cofibrations in $\spec$ are
small with respect to the whole category, see [HSS].
\end{proof}

The next lemma is used in proving Theorem~\ref{R-alg}.  It shows, 
for any commutative monoid $R$, that the category of 
$R$-modules has properties similar to the underlying category of $S$-modules. 

\begin{lemma}\label{R-mod+}
Let $R$ be a commutative ring in $\spec$.  Then the model category
structure on $R$-modules given above is a monoidal model category which 
satisfies the monoid axiom.   
\end{lemma} 

\begin{proof}
This follows from [SS, 4.1 (2)].
\end{proof}

\begin{theorem}\label{R-alg}
Let $R$ be a commutative monoid in $\spec$.  Then the category of $R$-algebras
is a cofibrantly generated model category with weak equivalences and
fibrations given by the underlying stable model category structure of $\spec$.  
The generating
cofibrations and trivial cofibrations are given by applying the free monoid 
functor to the generating (trivial) cofibrations of $R$-modules.
Moreover, if $f:A \xrightarrow{} B$ is a cofibration of $R$-algebras with
$A$ cofibrant as an $R$-module then $f$ is also a cofibration in the underlying
category of $R$-modules.  In particular, this shows that any cofibrant 
$R$-algebra is also cofibrant as an $R$-module. 	
\end{theorem}

\begin{proof}
This follows from Theorem 4.1 (3) in [SS].  The facts about cofibrant objects
and cofibrations of $R$-algebras follow from [SS, 4.1 (3)] because
$S$ is cofibrant in $\spec$.  
\end{proof}

Since symmetric ring spectra are exactly the $S$-algebras, the following
is a corollary of Theorem~\ref{R-alg}.  

\begin{corollary}
The category of symmetric ring spectra, $\Salg$, is a cofibrantly generated
model category with weak equivalences and fibrations the underlying stable 
equivalences and stable fibrations of $\spec$.  
\end{corollary}

The following lemma is needed to apply the comparison theorems of [SS, 4.3-4]
which show that a weak equivalence of symmetric ring spectra induces
an equivalence on the respective homotopy theories of modules and algebras.
This lemma is also needed for section~\ref{thh.1}. 

\begin{lemma}[HSS]\label{R-cof.^stable=stable}
Let $R$ be a symmetric ring spectrum and $M$ a cofibrant $R$-module.
Then $M \smashR -$ takes level equivalences of $R$-modules to level 
equivalences in $\spec$ and it takes stable equivalences of 
$R$-modules to stable equivalences in $\spec$.
\end{lemma}

The following two comparison theorems
follow from Lemma~\ref{R-cof.^stable=stable}, [SS, 4.3-4], and
the fact that $S$ is cofibrant in $\spec$.

\begin{theorem}If $A \xrightarrow{\sim} B$ is a map of
symmetric ring spectra which is an underlying stable equivalence, 
then the total derived functors
of restriction and extension of scalars induce equivalences of homotopy
theories
\[ \mbox{\em Ho}\,(A\mbox{\em -mod}) \ \iso \ \mbox{\em
Ho}\,(B\mbox{\em -mod}) \ . \]
\end{theorem}

\begin{theorem}\label{compare-algebras} 
If $A \xrightarrow{\sim} B$ is a map of commutative ring spectra which is
an underlying stable equivalence, then the total derived functors 
of restriction 
and extension of scalars induce equivalences of homotopy theories
 \[ \mbox{\em Ho}\,(A\mbox{\em -alg}) \ \iso \ \mbox{\em Ho}\,(B\mbox{\em -alg})
 \ . \]
\end{theorem}

The following two lemmas from [HSS] are used to verify some of the properties 
of the smash product that we have mentioned here.  They are also needed in
section~\ref{thh.1}.  

\begin{lemma}\label{cof^level=level}
Let $X \xrightarrow{} Y$ be a level equivalence.  Then $A \smashs  X
\xrightarrow{} A \smashs Y$ is a level equivalence for any cofibrant
spectrum $A$.
\end{lemma}

\begin{lemma}\label{cof^stable=stable}
Let $X \xrightarrow{} Y$ be a stable equivalence.  Then $A \smashs  X
\xrightarrow{} A \smashs Y$ is a stable equivalence for any cofibrant
spectrum $A$.
\end{lemma}

\subsection{Homotopy colimits}\label{hocolim}
In this section we list some of the properties of the homotopy colimit
functor for symmetric spectra which are used in the latter parts of this paper.
The most important property is that the homotopy colimit in symmetric
spectra can be defined by using the homotopy colimit of spaces at each
level, see Definition~\ref{hocolim levels}.  This is useful not only
because the homotopy colimit of spaces is well understood, but also
to show that the homotopy colimit preserves level equivalences of
symmetric spectra, see Proposition~\ref{level comparison}.
We use the basic construction of the homotopy colimit for spaces 
from [BK].

\begin{definition} \label{hocolim levels}
Let $B$ be a small category and $F: B \xrightarrow{}
\spectra$ a diagram of symmetric spectra.  Let $F_l$ denote the diagram
of spaces at level $l$.  Then
$$ (\hocolim_{\spectra}^B F)_l = \hocolim_{\ssets}^B F_l.$$
\end{definition}

This definition makes sense because any stable cofibration is a level
cofibration and colimits in $\spec$ are created on each level.  Also,
we show that this homotopy colimit has the usual properties of a
homotopy colimit.  Namely, a  map between diagrams which is objectwise
a level equivalence, a \p-isomorphism, or a stable equivalence
induces the same type of equivalence on the homotopy colimit.
The next two propositions consider the first two cases.  The case of
stable cofibrations could be proved by generalizing [BK, XII 4.2] to
arbitrary model categories.  Instead, here we use the detection functor
developed
in section~\ref{D} to verify this property in Lemma~\ref{stable hocolim}.

\begin{proposition}\label{level comparison}
Let $F,G :B \xrightarrow{} \spectra$ be two diagrams of symmetric spectra with
a natural transformation $\eta: F \xrightarrow{}G$ between them.  If
$\eta(b):F(b) \xrightarrow{} G(b)$ is a level equivalence at each object
$b \in B$, then $\hocolim^B F \xrightarrow{} \hocolim^B G$ is a level
equivalence.
\end{proposition}

\begin{proof}
Using Proposition~\ref{hocolim levels} this statement reduces to asking
that the homotopy colimit preserve objectwise weak equivalences of simplicial
sets.  This is the dual of [BK, XI 5.6].  Cofibrancy conditions are not
required here since any space (i.e., simplicial set) is cofibrant.
\end{proof}

We also need to know that the homotopy colimit of an objectwise
\p-isomorphism is a \p-isomorphism.
For this we form a spectral sequence for calculating the classical stable
homotopy groups
of the homotopy colimit.   Following [BK, XII 5] one can form a spectral
sequence for calculating any homology theory applied to the homotopy colimit
of spaces.
The spectral sequence is associated to the filtration of the homotopy colimit
given by the length of the sequence of maps in $B$.
So for $F\s B \to \sset$
this spectral sequence converges to $h_* \hocolim^B F$ and has $E^2$-term
$$ E^2_{s,t}= \colim^s_B  (h_t F).$$

We use the following lemma to go from the homology theory $\pi_*^s$ defined on
spaces by  $\pi_*^s K = $\p$F_0K$  to one on $\spec$.

\begin{lemma}
For $X$ a symmetric spectrum, \p$X=\colim_n \pi_*^s X_n$
\end{lemma}

\begin{proof} Consider the lattice of spaces 
$\Omega^i(\Omega^jL'\Sigma^j X_i)$
indexed over $(i, j) \in \N \times \N$ with maps for fixed $j$
using the adjoint structure maps of $\Omega^jL' \Sigma^j X$ or
for fixed $i$ using $\Omega^i\Omega^j$ applied to the adjoint structure
maps of $L'F_0 X_i$.  Applying homotopy and taking colimits in the two different
directions finishes the proof.  In one direction, one gets $\colim_j
\pi_* \Omega^j L'\Sigma^j X$, but each of these terms and hence the colimit is 
isomorphic to $\pi_* X$.   In the other direction, one has the $\colim_i
\pi_*^s X_i$.  
\end{proof}

So applying the homology theory \p, the above spectral sequence calculates
\p\  of
each level of the homotopy colimit.  Since \p$X=\colim_n \pi_*^s X_n$
and a sequential colimit of spectral sequences is a spectral sequence,  taking
the colimit of these level spectral sequences produces a spectral
sequence.

\begin{proposition}\label{spec.seq} For $F\s B \to \spec$, there is a 
spectral sequence converging to \p $ \hocolim^B_{\spec} F$ with $E^2$-term
$$ E^2_{s,t}= \colim^s_B (\pi_t F).$$
\end{proposition}

Using this spectral sequence we can show that homotopy colimits preserve
objectwise \p-isomorphisms.

\begin{proposition}\label{stable comparison}
Let $F,G :B \xrightarrow{} \spectra$ be two diagrams of symmetric spectra with
a natural transformation $\eta: F \xrightarrow{}G$ between them.  If
$\eta(b):F(b) \xrightarrow{} G(b)$ induces a \p-isomorphism
at each object
$b \in B$ then $\hocolim^B F \xrightarrow{} \hocolim^B G$ induces
a \p-isomorphism.
\end{proposition}

\begin{proof}
Since $\eta$ induces a \p-isomorphism between the two diagrams in question,
it induces an $E^2$-isomorphism.  Thus it induces an isomorphism on
the $E^\infty$-term and hence a \p-isomorphism on the homotopy colimits.
\end{proof}

For section~\ref{D} we also need the following two cofinality
results which are from [BK, XI 9.2].

Given a functor between two small categories
$f: A \xrightarrow{} B$ one
has a natural map $\hocolim^A f^*F \xrightarrow{} \hocolim^B F$.  A functor
$f:A \xrightarrow{}  B$ is called {\em terminal} or {\em right cofinal}, if
for every object $b \in B$ the under category $(b\-\downarrow\- f)$ is
contractible see [BK, XI 9].

\begin{proposition}\label{terminal}
Let $f:A \xrightarrow{} B$ be a functor which is terminal.  Then
for any functor $F:B \xrightarrow{} \spectra$,
$\hocolim^A f^*F \xrightarrow{} \hocolim^B F$ is a level equivalence.
\end{proposition}

\begin{proof}[Proof.]
The dual of [BK, XI 9.2] states this property for objectwise
cofibrant diagrams of spaces.  Applying this on each level and using the
fact that spaces are cofibrant proves this statement.
\end{proof}

In section~\ref{D}, we consider diagrams over the skeleton of the category
of finite sets and injections, $I$, with objects ${\bf n}$.  
Let $I_{m}$ denote the full
subcategory of $I$ whose objects are ${\bf n}$ where $n$ is
greater than or equal to $m$.  The following
lemma states the cofinality information relating these categories.

\begin{lemma}\label{cofinal} Let $F\s I \to \spec$ be a diagram of spectra.
The inclusion $u_m:I_{m} \xrightarrow{} I$ is terminal, hence
$\hocolim^{I_{m}}u_m^* F \to \hocolim^{I} F$ is a level equivalence.
\end{lemma}

\begin{proof}
Consider the functor $- + m:
I \xrightarrow{} I_{m}$ which induces a functor on the under categories.
There is a natural transformation from the identity functor to both
$u_m \circ (- +m)$ and $(- +m ) \circ u_m$.
Hence the under categories are each homotopy equivalent to $(i\downarrow I)$.
But $(i\downarrow I)$ is contractible because it has an
initial object $1: i\xrightarrow{} i$.  The homotopy colimit statements
follow from Proposition~\ref{terminal}.
\end{proof}

Using these cofinality results we can prove the following proposition.

\begin{proposition}\label{connectivity}
Let $F,G \s I \xrightarrow{} \ssets$ be two diagrams of spaces with
a natural transformation $\eta: F \xrightarrow{}G$ between them.  Assume that
$\eta({\bf n})\s F({\bf n}) \xrightarrow{} G({\bf n})$ is a $\lambda(n)$
connected map, where $\lambda(n) \leq \lambda(n+1)$ and $lim_n \lambda(n)$ is
infinite. Then
$\hocolim^I F \xrightarrow{} \hocolim^I G$ is a weak equivalence.
\end{proposition}

\begin{proof}
We show that the map is an $N$-equivalence for every $N>0.$  Choose an $n$
such that $\lambda(n) > N$.  Then for every object ${\bf m}$ in $I'_{n}$
the map $\eta \s F({\bf m}) \to G({\bf n})$ is an $N$-equivalence, and
so we can conclude that $\eta\s \hocolim^{I_{n}} u_n^* F
\to \hocolim^{I_{n}} u_n^* G$ is an $N$-equivalence.  The proposition
follows by Lemma~\ref{cofinal}.
\end{proof}

We also need the following proposition which shows that the homotopy colimit of
a diagram of level equivalences over $I$ is level equivalent to its value
on ${\bf 0}$.

\begin{proposition}\label{initial}
Let $F\s I \to \spec$ be a diagram of spectra.  Assume that
for each morphism $f$ in $I$, $F(f)$ is a level equivalence.  Then
the inclusion $F({\bf 0}) \to \hocolim^I F$ is a level equivalence.
\end{proposition}

\begin{proof}
Consider the constant functor $C\s I \to \spec$ with constant value
$F({\bf 0})$.  Then at each object the map $C({\bf n})=F({\bf 0}) \to
F({\bf n})$ induced by the unique map ${\bf 0} \to {\bf n}$ is a level
equivalence.  Hence, by Proposition~\ref{level comparison},
it induces a level equivalence on
the homotopy colimits, $F({\bf 0}) \to \hocolim^I F$.
\end{proof}

Finally, we need the following proposition due to Jeff Smith, [S].
Let $T$ be the category with objects ${\bf{n}}=\{1, \dots, n\}$ and morphisms 
the standard inclusions.  Homotopy colimits over $T$ are weakly equivalent
to telescopes.  
Let $\omega$ be the ordered set of natural numbers 
and $I_{\omega}$ be the category whose objects are the finite sets ${\bf n}$
and the set $\omega$ and whose morphisms are inclusions.  
Let $L_h F\s I_{\omega} \to 
\ssets$ be the left homotopy Kan extension of $F\s I \to \ssets$  along
the inclusion of categories $i \s I \to I_{\omega}$.  

\begin{proposition}\label{smith}  
Let $M$ be the monoid of injective maps $i\s \omega \to \omega$ under 
composition. Given any functor $F \s I \to \ssets$, then 
\begin{enumerate}
\item $\hocolim^I F$ is weakly equivalent to $(L_hF(\omega))_{hM}$ where
$(-)_{hM}$ is the homotopy orbits with respect to the action of $M$, and
\item $L_h F(\omega)$ is weakly equivalent to $\hocolim^T F$.
\end{enumerate}     
\end{proposition}

\begin{proof} 
For the convenience of the reader we sketch Smith's proof of this proposition.
Since $L_h F$ is the homotopy Kan extension, 
$\hocolim^I F \hoeq \hocolim^{I_\omega} L_h F$.
Next, consider the full subcategory, $A$ of $I_{\omega}$ with just one object,
$\omega$.  Since the inclusion of $A$ in $I_{\omega}$ is terminal,   
$\hocolim^{I_{\omega}} L_hF$ is weakly equivalent to $\hocolim^{A} L_hF$. 
Since $\Hom_{A}(\omega, \omega)=M$, 
$\hocolim^{A} L_hF$ is the homotopy orbit space $(L_h F(\omega))_{hM}$.  

For the second statement, 
$L_hF(\omega)= \hocolim^{i(n) \to \omega \in (i\-\downarrow\-\omega)} F(n)$.  
Here $i$ is the inclusion $I \to I_{\omega}$.  The category
$T$ described above is equivalent to the category 
$(i\circ \alpha\-\downarrow\-\omega)$ for the inclusion $\alpha \s T \to I$.
This category $(i\circ \alpha\-\downarrow\-\omega)$ is
terminal in $(i\-\downarrow\-\omega)$, because every
under category has an initial object.  So by Proposition
\ref{terminal}, $L_hF(\omega)$ is weakly equivalent to $\hocolim^{T} F$. 
\end{proof}

\section{Detecting stable equivalences}\label{D}
In this section we introduce a functor, $D$, which detects stable 
equivalences in the sense that a map $X \xrightarrow{} Y$ is a stable 
equivalence
if and only if $DX \xrightarrow{} DY$ is a \p-isomorphism.
Of course the stable fibrant replacement functor, $L$, also has this property.
It even turns stable 
equivalences into level equivalences.  The drawback of $L$ is that its only 
description is via the small object
argument.  Hence it is difficult to say much about $L$ apart from its 
abstract properties.  The advantage of the functor $D$ is
that it has a more explicit definition.  In particular, there is a spectral
sequence for calculating the classical stable homotopy groups of $DX$, see
Proposition \ref{spec.seq}.   Moreover, these
groups are invariants of the stable equivalence type of $X$ because $D$ 
takes stable equivalences to \p-isomorphisms.   

In Section~\ref{thh.1} we see that $D$ fits into a sequence of functors used
to define $\THH$ in [B].  
We use the notation $D$ instead of $\THH_0$ because $D$ is defined on any
symmetric spectrum, not just on ring spectra. 

\subsection{Main statements and proofs}\label{D.1}
The detection functor $D$ is defined as a homotopy colimit over the diagram
category of the skeleton of finite sets and injections, $I$.  
Given a symmetric spectrum $X$, define a functor $\mathcal{D}_X\s I 
\xrightarrow{} \spectra$ whose value on the object ${\bf n}$ is 
$\Omega^n L' F_0X_n$.  Recall $L'$ is just a level fibrant replacement functor.
For a standard inclusion of a subset 
$\alpha \s {\bf n} \subset {\bf m}$ the map ${\mathcal{D}}_X(\alpha)$ is 
just $\Omega^n L'$ applied to the composition of  maps $F_0X_n \xrightarrow{} 
F_0 \Omega^{m-n}X_m
\xrightarrow{} \Omega^{m-n}F_0 X_m$ induced by the structure maps of $X$.
For an isomorphism, the action is given by the conjugation action on
the loop coordinates and on $X_n$. 
All morphisms in $I$ are compositions of isomorphisms and these standard 
inclusions.

\begin{definition}\label{def D} The \emph{detection functor} 
$D:\spectra \xrightarrow{} \spectra$
is defined by $$ DX= \hocolim^{I}_{\spec}  \mathcal {D}_X.$$  
\end{definition}

As defined in definition~\ref{hocolim levels}, the homotopy colimit of 
symmetric spectra
is given by a level homotopy colimit of spaces.  
Hence $$(DX)_n=\hocolim^{k\in I}_{\ssets} \Omega^k L' \Sigma^n X_k.$$

As mentioned above, the main reason for considering $D$ is that it detects 
stable equivalences.  This is stated in the next theorem.

\begin{theorem}\label{detection-thm} 
The following are equivalent.
\begin{enumerate}
\item $X \xrightarrow{} Y$ is a stable equivalence.
\item $DX \xrightarrow{} DY$ induces a \p-isomorphism.
\item $D^2X \xrightarrow{} D^2Y$ is a level equivalence.
\item $QDX \xrightarrow{} QDY$ is a level equivalence.
\end{enumerate}
\end{theorem}

\begin{remark} Notice that one can apply the forgetful functor 
$U\s \spec \to \BF$ after applying $D$.  In that case, this theorem 
says that although the usual
forgetful functor does not detect and preserve stable  
equivalences, the composition of this detection functor with the forgetful
functor does detect and preserve weak equivalences.  
Note also that although the classical stable homotopy groups are not 
invariants of stable equivalence types 
in symmetric spectra this theorem shows that  after applying $D$ 
the classical stable homotopy groups are invariants.
\end{remark}

\begin{remark} \label{top-spaces}
We could also consider symmetric spectra over topological spaces
instead of simplicial sets here, see [HSS].  In that case, Theorem~\ref{detection-thm}
and all of the statements leading up to it in this section and in section
\ref{hocolim} which do not involve the functor $Q$ hold when the objects
involved are levelwise non-degenerately based spaces.  Hence,
$D$ also detects stable equivalences between symmetric spectra based on
topological spaces.  More precisely, let $c$ be a cofibrant replacement functor
of spaces applied levelwise, then $X \to Y$ is a stable equivalence if and
only if $DcX \to DcY$ is a $\pi_*$-isomorphism.

The only fact that is needed to modify all of these statements for topological
spaces is that homotopy colimits of non-degenerately
based spaces are invariant under weak homotopy equivalences.
For the statements involving $Q$ one needs to consider stably cofibrant
symmetric spectra because these statements require that homotopy groups
commute with directed colimits.  But these statements are separate from the
statements involving $D$.
\end{remark}

Theorem~\ref{detection-thm} considers the properties of $D$ with respect
to morphisms.  The following theorem considers the properties of $D$ on
objects.

\begin{theorem}\label{thm-objects}
Let $X$ be a symmetric spectrum.
\begin{enumerate}
\item $DX$ is semistable.
\item If $X$ is semistable, then the level fibrant replacement
of $DX$, $L'DX$, is an $\Omega$-spectrum.
\end{enumerate}
\end{theorem}

Since stable equivalences between semistable spectra are \p-isomorphisms
and between $\Omega$-spectra are level equivalences,  Theorem \ref{thm-objects}
shows that the second and third statements of Theorem~\ref{detection-thm} 
really just say that $D$ and $D^2$ preserve and detect stable equivalences.

Theorem~\ref{detection-thm} shows that the classical stable homotopy groups
of $DX$ are a stable equivalence invariant.  In the next theorem we show
that they are in fact the derived classical stable homotopy groups, {\it i.e.,}
they are isomorphic to \p$LX$. 

\begin{theorem}\label{thm-derived}
Let $X$ be a symmetric spectrum.
\begin{enumerate}
\item There is a natural zig-zag of functors
inducing \p-isomorphisms between $LX$ and $DX$.  
\item There are natural zig-zags of functors 
inducing level equivalences between $LX$, $D^2X$, and $QDX$.
\end{enumerate}
\end{theorem}

This theorem shows that the fibrant replacement functor is determined up to 
\p-isomorphism by $D$ or up to level equivalence by $D^2$ or $QD$.  
The spectral sequence for calculating the classical stable homotopy groups 
of $DX$, Proposition \ref{spec.seq}, calculates the derived 
stable homotopy groups \p$DX \cong$ \p$LX$.

\begin{corollary}\label{cor-sf}
For $X$ any semistable spectrum, $X$ and $DX$ are \p-isomorphic.
\end{corollary}  

\begin{remark}\label{sf-remark}
This corollary is a special case of [B, 1.6] where the convergence
and connectivity conditions are replaced by the semistable condition. 
By Proposition~\ref{sf-ex} we recover a statement with convergence
conditions but no connectivity conditions.
\end{remark}  

The proofs of Theorems~\ref{detection-thm} and \ref{thm-objects} 
use the following properties of the functor $D$.

\begin{proposition}\label{detection-prop}
Let $f:X \xrightarrow{} Y$ be a map of symmetric spectra.
\begin{enumerate}
\item If $f$ is a stable equivalence then $Df$ is a \p-isomorphism. 
\item If $f$ is a \p-isomorphism then $Df$ is a level equivalence. 
\item For any semistable spectrum $X$, there is a natural zig-zag of
functors inducing level equivalences between $LX$ and $DX$.
\end{enumerate}
\end{proposition}

We assume Proposition~\ref{detection-prop}
to prove Theorems~\ref{detection-thm}, \ref{thm-objects}, and \ref{thm-derived}.
The proof of Proposition~\ref{detection-prop} 
is technical, so it is delayed until the next subsection.

\begin{proof}[Proof of Theorem~\ref{thm-derived}]  
By Proposition~\ref{detection-prop} (3) applied to $LX$
there is a zig-zag of level equivalences between $LLX$ and $DLX$.  By
Proposition~\ref{detection-prop} (1) since $X \to LX$ is a stable equivalence
$DX \to DLX$ is a \p-isomorphism.   Hence, putting these equivalences
together with the fact that $LLX$ is level equivalent to $LX$, we get
a zig-zag of \p-isomorphisms between $LX$ and $DX$.

Applying $D$ to the zig-zag of \p-isomorphisms between $LX$ and $DX$
shows that $DLX$ and $D^2X$ are level equivalent by 
Proposition~\ref{detection-prop} (2).  Combining this with the zig-zag
of level equivalences between $LX$ and $DLX$ produces the level equivalence of
$LX$ and $D^2X$. The equivalences for $QDX$ are similar.
\end{proof}

\begin{proof}[Proof of Theorem \ref{thm-objects}]
By Theorem \ref{thm-derived} $DX$ is \p-isomorphic to $LX$.  $LX$ is an
$\Omega$-spectrum, hence it is semistable.  So by Proposition \ref{sf-p-eq},
$DX$ is semistable.

For $X$ semistable, Proposition \ref{detection-prop} shows that
$DX$ is level equivalent to $LX$, an $\Omega$-spectrum.  Hence $L'DX$ is
an $\Omega$-spectrum. 
\end{proof}

\begin{proof}[Proof of Theorem \ref{detection-thm}]  
Proposition \ref{detection-prop} shows that (1) implies (2) and (2) implies (3).
A map $f$ is a \p-isomorphism 
if and only if $Qf$ is a level equivalence.  Hence the second 
and fourth statements are also equivalent.  

By Theorem \ref{thm-derived} part 2, $LX$ and $D^2X$ are naturally level 
equivalent.
Hence if $D^2X \to D^2Y$ is a level equivalence then so is $LX \to LY$.
But this is equivalent to $X \to Y$ being a stable equivalence.
\end{proof}

\subsection{Proof of Proposition \ref{detection-prop} }
As mentioned above the proof of 
Proposition \ref{detection-prop} is more technical.  In this subsection
we first prove the second part of Proposition \ref{detection-prop}.  Using
this we prove the third part.  Then, for the first part of Proposition 
\ref{detection-prop} we state and prove several lemmas which together finish
the proof.  Throughout this section we use several of the properties of the 
homotopy colimit developed in section \ref{hocolim}.

For the proof of the second part of Proposition \ref{detection-prop} we
use Proposition \ref{smith}, due to Jeff Smith.  

\begin{proof}[Proof of Proposition \ref{detection-prop} Part 2]
We apply Lemma \ref{smith} to each level of $D$.  Consider
the $0$th level first.  If $f$ is a \p-isomorphism then 
$\hocolim^T \Omega^n L'f_n$  is a weak equivalence,
since \p$X= $\p$\hocolim^T\Omega^nL'X_n$.
Since taking homotopy orbits preserves weak equivalences this shows
that  the $0$th level of $DX \to DY$ is a weak equivalence, {\it i.e.,} 
$\hocolim^I \Omega^n L'f_n$ is a weak equivalence.

The $k$th level of $DX$ is the $0$th level of $D \Sigma^k X$. 
Since $\Sigma^k f$ is a \p-isomorphism if $f$ is, this shows that each
level is a weak equivalence.
\end{proof}

For the third part of Proposition \ref{detection-prop} we need the following
functor.  Recall that $(sh^n X)_k = X _{n+k}$.    

\begin{definition}\label{def M} Define $MX = \hocolim^{I} 
\Omega^n L'sh_n X$.  
\end{definition} 
 
\begin{proof}[Proof of Proposition \ref{detection-prop} Part 3]
First we develop the transformations which play
a part in the zig-zag mentioned in the proposition. 
The inclusion of the object ${\bf 0}$ in $I$
induces a natural map $X \xrightarrow{} MX$.  There is also a natural
transformation of functors $D \xrightarrow{} M$.  
The structure maps on $X$ induce a natural map of symmetric spectra 
$F_0X_n \xrightarrow{}
sh_n X$.   Applying $\Omega^n L'$ this map induces a map of diagrams over $I$,  
and hence a natural map of homotopy colimits.
So there is a natural zig-zag $X \xrightarrow{} MX \xleftarrow{} DX$.
The zig-zag mentioned in the proposition is this zig-zag applied to $LX$ along
with the natural map $DX \to DLX$.  

For semistable $X$, the map $X \to LX$ is a \p-isomorphism.  So $DX \to DLX$
is a level equivalence by Proposition \ref{detection-prop} part 2.  So
we only need to show that if $X$ is an $\Omega$-spectrum, then both of the maps 
$X \to MX \xleftarrow{} DX$ are level equivalences.  

First we show that the map $X \xrightarrow{} MX$ is a level equivalence for
any $\Omega$-spectrum $X$.  By definition an $\Omega$-spectrum
is a  level fibrant spectrum such that $X \xrightarrow{} \Omega \sh_1 X$ is a
level equivalence.  Using this and the fact that both shift and $\Omega$ 
preserve level equivalences (on level fibrant spectra), one can show 
that each of the maps in the
diagram over $I$ used to define $MX$ is a level equivalence.  
By Proposition \ref{initial} this implies that $X \xrightarrow{} MX$ is a level 
equivalence.       

To show that $DX \xrightarrow{} MX$ is a level equivalence for any 
$\Omega$-spectrum $X$, we need to consider connective covers.
Given a  level fibrant spectrum $X$ define its $k$th connective
cover, $C_kX$,  as the homotopy fibre of the map from $X$ to its $k$th
Postnikov stage $P_kX$.  
The $k$th Postnikov functor is the localization
functor given by localizing with respect to the set of maps 
$\{F_n \partial\Delta[m+n+k+2] \to F_n \Delta[m+n+k+2] 
\s m, n \geq 0\}$. 
At level $n$, this functor is weakly equivalent to the $(n+k)$th Postnikov 
functor on spaces  which is given by localization with respect to the set of 
maps $\{\partial\Delta[m+n+k+2] \to \Delta[m+n+k+2]\s m\geq 0\}.$  See also [F2].  
Then $(C_kT)_n$ is $n+k$ connected and
$\pi_i (C_kT)_n \to \pi_i T_n$ is an isomorphism for $i > n+k$.     
Note that any level fibrant spectrum is level equivalent to the homotopy 
colimit over its connective covers.  As $-k$ decreases, 
the homotopy type of each level of $C_{-k}X$ eventually becomes constant.  
So $\hocolim_{k} C_{-k}X \to X$ is a level equivalence.

Because $\Omega^m$, $L'$ and $F_0$ commute up to level equivalence with 
directed homotopy colimits and homotopy 
colimits commute, $\hocolim_n DC_{-n}X$ is level equivalent to $DX$.  The shift
functor also commutes with homotopy colimits so $\hocolim_n MC_{-n}X$ is
level equivalent to $MX$.  
So we first show that for each $n$, $DC_{n}X$ and $MC_{n}X$ are 
level equivalent.  

In the diagrams creating these homotopy colimits, consider level $l$ at
the object $m \in I$.  The map in question is $\Omega^m L' (S^l \wedge 
C_{n}X_m) \xrightarrow{} \Omega^m L' C_nX_{m +l} $.  In general the map
$\Omega^m L' (\Sigma^l \Omega^l L' Y) \xrightarrow{} \Omega^m L' Y$ is 
$2N -l - m + 1$ connected when $Y$ is $N$ connected. 
Hence for $Y=C_nX_{m + l}$ the map in question is $2n+m+l+1$ connected.  

Using Proposition \ref{connectivity}, we see that this connectivity implies 
that $(DC_nX)_k \xrightarrow{} (MC_nX)_k$ is a weak equivalence.   
Homotopy commutes with directed colimits, so taking colimits over $n$ on
both sides we get a weak equivalence $DX_k \xrightarrow{} MX_k$.  So
$DX \xrightarrow{} MX$ is a level equivalence.  This is what we needed to
finish the third part of Proposition \ref{detection-prop}.  
Note also that since $MX$ is 
an $\Omega$-spectrum this shows that the level fibrant replacement of $DX$ is
also an $\Omega$-spectrum.  
\end{proof}

The proof of Proposition \ref{detection-prop} part 1 breaks up into several 
parts.  For the case of stable trivial cofibrations we split the problem
into showing that $D$ of any generating stable trivial cofibration is a 
\p-isomorphism and that $D$ behaves well with respect to push 
outs, {\it i.e.,} that the following two lemmas hold.

\begin{lemma}\label{D-J}
Let $j: A \xrightarrow{} B$ be a generating stable trivial cofibration.  Then
$Dj: DA \xrightarrow{} DB$ is a \p-isomorphism.
\end{lemma}

\begin{lemma}\label{D-p.o.}
If 
\[
\begin{CD}
A @>>> X\\
@VVV   @VVV\\
B @>>> Y
\end {CD}
\]
is a pushout square with $A \xrightarrow{} B$ a cofibration, then
\[
\begin{CD}
DA @>>> DX\\
@VVV   @VVV\\
DB @>>> DY
\end {CD}
\]
is a homotopy pushout square.  I.e, if $P$ is the homotopy colimit 
of $DB \xleftarrow{} DA \xrightarrow{} DX$, then $P \xrightarrow{} DY$ is 
a stable equivalence.
In fact, $P \xrightarrow{} DY$ is a \p-isomorphism.
\end{lemma}

Combining this lemma with the next shows that if $DA \xrightarrow{} DB$ is
a \p-isomorphism then $DX \xrightarrow{} DY$ is also a \p-isomorphism.

\begin{lemma}\label{stable proper}
Let \[
\begin{CD}
A @>>> X\\
@VVV   @VVV\\
B @>>> Y
\end {CD}
\]
be a square in $\spec$ with $Y$ \p-isomorphic to the homotopy pushout. 
Assume $A \xrightarrow{} B$ is a \p-isomorphism.  
Then $X \xrightarrow{} Y$ is a \p-isomorphism.
\end{lemma}

For a proper model category this is a standard fact, that the homotopy
pushout of a weak equivalence is a weak equivalence.  But no model category
on symmetric spectra has been written down with weak equivalences the
\p-isomorphisms.

\begin{proof}[Proof of Proposition \ref{detection-prop} Part 1]
Assuming Lemmas \ref{D-J}, \ref{D-p.o.}, and \ref{stable proper},
we can finish this proof.
First note that since any stable equivalence can be factored as a stable trivial
cofibration followed by a level trivial fibration, we only need to show that 
$D$ takes both stable 
trivial cofibrations and level equivalences to \p-isomorphisms. 

A level equivalence induces a level equivalence at each object in the diagram
for defining $D$ and homotopy colimits preserve level equivalences,
by Proposition\ref{level comparison}.   Hence, $D$ of a level equivalence is 
a level equivalence, and thus a \p-isomorphism. 
 
Any stable trivial cofibration is a retract of a directed colimit of
pushouts of maps in $J$.
Since retracts and directed colimits preserve 
\p-isomorphisms we only need to consider pushouts of generating stable trivial
cofibrations.  By Lemma \ref{D-J}, $D$ of a generating trivial cofibration 
is a \p-isomorphism.  Hence, by Lemmas \ref{D-p.o.} and \ref{stable proper},
$D$ of any map formed by a pushout of a generating stable trivial
cofibration is a \p-isomorphism.  
\end{proof} 

\begin{proof}[Proof of Lemma \ref{stable proper}]
Factor the map $A \xrightarrow{} X$ as a stable cofibration  
followed by a level trivial fibration $A \to Z \xrightarrow{} X$.  Then form 
the pushout square as follows,    
\[
\begin{CD}
A @>>> Z \\
@VVV    @VVV\\
B @>>> P. 
\end{CD}
\]
Since the top map is a level cofibration, $P$ is the homotopy pushout of this 
square.  Since $A \xrightarrow{} B$ is a \p-isomorphism, 
$Z \xrightarrow{} P$ is a \p-isomorphism because \p\ is a homology theory.  
Since $Z \xrightarrow{} X$ is a level equivalence, to see that 
$X \xrightarrow{}Y$ is a \p-isomorphism it 
is enough to know that $P \xrightarrow{} Y$ is a \p-isomorphism.
But this is assumed as part of the hypotheses. 
\end{proof}

Now we proceed with the proof of Lemma \ref{D-p.o.}.

\begin{proof}[Proof of Lemma \ref{D-p.o.}]
To see that $P \xrightarrow{} DY$ is a \p-isomorphism,
we use the fact that homotopy colimits commute.  $P$ is the homotopy colimit 
of $DB
\xleftarrow{} DA \xrightarrow{} DX$, 
so $P$ is level equivalent to the homotopy colimit
over $I$ of the homotopy pushout of the squares at each object in $I$.   
In other words,  let $P^n$ be the homotopy pushout at the object 
${\bf n} \in I$. Then $P$ is level equivalent to $\hocolim^I P^n$.  

Proposition \ref{stable comparison} shows that a map of diagrams 
which is a \p-isomorphism at each object induces a 
\p-isomorphism on the homotopy colimits. 
Hence, it is enough to show that $P^n \xrightarrow{}
\Omega^n L' F_0(Y_n)$ is a \p-isomorphism for each $n$.  

Since cofibrations induce level cofibrations and $F_0$ preserves 
cofibrations and pushouts,  $F_0$ applied to each level of the pushout square
in the lemma is a homotopy pushout square. 
Since $X \to L'X$ is a level equivalence it preserves homotopy pushout
squares up to level equivalence.  Since $\Omega^n$ only shifts \p\  by $n$,
it preserves homotopy pushouts up to \p-isomorphism.  Hence $P^n
\to \Omega^n L' F_0(Y_n)$ is a \p-isomorphism.
\end{proof}

We are left with proving Lemma \ref{D-J}.  First we prove the following lemma 
which identifies the stable homotopy type of $DF_m(K)$.  

\begin{lemma}\label{stable type}
There is a $2l -m-1$ connected map $\psi_l \s
\Omega^mL'(S^l \sm K) \to (DF_mK)_l$. 
These maps fit together to give a map of symmetric spectra $\psi \s
\Omega^m L'F_0 K \to DF_mK$ which is a \p-isomorphism. 
\end{lemma}

To prove this lemma we define another functor on the category $I$.

\begin{definition}
Define $\F_mK \s I \to \sset$ by $(\F_mK)(n)=\hom_I(m,n)_+ \sm K$.  
\end{definition}

$\F_m(-)$ is left adjoint to the functor from $I$-diagrams over
$\sset$ to $\ssets$ which evaluates the diagram at $m \in I$.  Hence 
a natural transformation from $\F_mK$ into any diagram over $I$ is determined
by a map from K to the diagram evaluated at $m$.

\begin{proof} [Proof of Lemma \ref{stable type}] 
Let $\D^l_{F_mK}\s I \to \sset$ be the functor given by the $l$th 
level of the functor $\D_{F_mK}$.  Then there is a map  
$\phi_l\s \F_m\Omega^mL'(S^l \sm K) \to {\D}^l_{F_mK}$ determined by 
the inclusion 
of the wedge summand corresponding to the identity map,  $\Omega^mL'(S^l \sm K) 
\to \Omega^mL'(S^l \sm \hom_I(m,m)_+ \sm K).$  Because the
homotopy colimit of a free diagram is weakly equivalent to the colimit, 
see [F1], the homotopy colimit of this map is the map $\psi_l \s
\Omega^mL'(S^l \sm K) \to (DF_mK)_l$ mentioned in the lemma.

We show that the map of diagrams
is $2l-m-1$ connected at each spot.  At each $n\in I$, $\phi_l(n)$, factors 
into two maps as follows, 
\[\hom_I(m,n)_+ \sm \Omega^mL'(S^l \sm K) \to \Omega^mL'(\hom_I(m,n)_+ \sm S^l \sm K)  
\]
\[
\to 
\Omega^m \Omega^{n-m}L' \Sigma^{n-m}(\hom_I(m,n)_+ \sm S^l \sm K).\]
The first map is $2l -m - 1$ connected by the application of the 
Blakers-Massey theorem which shows that a wedge of loop spaces, $\Omega X \vee
\Omega Y$,  is equivalent in the stable range to the loop of the wedge, 
$\Omega(X \vee Y)$.  
The second map is $2l -m -1$ connected by the Freudenthal suspension theorem,
which for simplicial sets concerns the map $X \to \Omega L' \Sigma X$.
Hence the map at each spot in the diagram, $\phi_l(n)$ and thus the map of 
homotopy colimits, $\psi_l$ is $2l -m-1$ connected. 

To see that these levels fit together, note that we can prolong $\F_m$ to
a functor from symmetric spectra to $I$-diagrams of symmetric spectra.
Then there is a map $\phi\s 
\F_m(\Omega^m L'F_0K) \to \D_{F_mK}$ which on level $l$
is given by the map, $\phi_l$, above.  
Hence, taking homotopy colimits, this induces
a map $\psi\s \Omega^m L'F_0K \to DF_mK$ which is a \p-isomorphism. 
\end{proof}

\begin{proof}[Proof of Lemma \ref{D-J}] 
We must show that $D$ of a generating trivial 
cofibration is a \p-isomorphism.  First recall that the generating
trivial cofibrations are the maps $P(c, F_mg_r): P_{m,r} \xrightarrow{}
C \otime_S F_m(\Delta\K)$ where $P_{m,r}$ is the pushout below.  
$$\begin{CD}
F_1S^1 \otime_S F_m(\dot{\Delta}\K)  @>>> 
F_1S^1 \otime_S F_m({\Delta}\K) \\ 
@VVV     @VVV\\
C \otime_S F_m(\dot{\Delta}\K) @>>> P_{m,r}
\end{CD}$$
To show that
$D$ of $P(c, F_mg_r)$ is a \p-isomorphism it is only
necessary to show that $D$ of $c_K:F_1S^1 \otime_S F_mK \xrightarrow{}
C \otime_S F_mK$ is a \p-isomorphism for 
$K=\dot{\Delta}[r]_{+}$
or $\Delta[r]_{+}$.  This is enough because Lemma \ref{stable proper} shows that  
if $D(F_1S^1 \otime_S F_m(\dot{\Delta}[r]_+) ) \xrightarrow{} D(C \otime_S
F_m(\dot{\Delta}[r]_+))$ is a \p-isomorphism then the pushout 
$D(F_1S^1 \otime_S F_m(\Delta[r]_+) ) \xrightarrow{} D(P_{m,r})$
is also a \p-isomorphism.  If  
$D(F_1S^1 \otime_S F_m(\Delta[r]_+) ) \xrightarrow{} D(C \otime_S
F_m(\Delta[r]_+))$ is also a \p-isomorphism, this implies that
$D(P_{m,r}) \xrightarrow{} D(C \otime_S F_m(\Delta[r]_+) )$ is a \p-isomorphism.

Since $F_mK$ is cofibrant and $C \xrightarrow{} F_0S^0$ is a level 
equivalence, the map 
$C \otime_S  F_mK \xrightarrow{} F_0S^0 \otime_S F_mK$ is   
a level equivalence, by Lemma \ref{cof^level=level}.
As already noticed, $D$ takes level
equivalences to \p-isomorphisms so we can assume that $C$ is
replaced by $F_0S^0$ in $c_K$ for both values of $K$. 

Thus we only need to show that $Dc_K \s DF_{m+1}(S^1 \sm K) \to DF_mK$
is a \p-isomorphism.  
To see that this map induces an isomorphism on the
stable homotopy groups, note that $F_{m+1}(S^1 \sm K) \to F_mK$ is
induced by $\hom_I(m+1,n) \to \hom_I(m,n)$ which in turn is induced
by the inclusion of $m$ in $m+1$.  
Now consider homotopy applied to the map of diagrams, $\D c_K$.  Using the
\p-isomorphisms from Lemma \ref{stable type} above,  this map is a map of
free diagrams, $\hom_I(m+1, -) \otimes \pi^s_{*+m+1}S^1 \sm K \to  
\hom_I(m, -) \otimes \pi^s_{*+m} K$.  This map induces an isomorphism
on the colimits and all of the higher $\colim^i$ vanish.  Hence, using
the spectral sequence for calculating the homotopy of homotopy colimits,
see section~\ref{hocolim}, 
$Dc_K$ is a \p-isomorphism.  One can also see this by considering the
associated map of free diagrams directly.  
\end{proof}

\section{Topological Hochschild Homology}\label{thh.1}
Let $k$ be a commutative symmetric ring spectrum.  Let $R$ be a $k$-algebra.
Define $R^e = R \smashk R^{op}$.  
Let $M$ be a $k$-symmetric $R$-bimodule, {\it i.e.,} an $R^e$-module.
With this set up we have two different definitions of topological Hochschild
homology, one using a derived tensor product definition, the other
mimicking the usual Hochschild complex.  In Theorem \ref{thh comparison} 
we see that these definitions construct stably equivalent $k$-modules.
Of course, since the smash product is only stably invariant for cofibrant
spectra, the case where $R$ is a cofibrant $k$-module is the only one of 
interest.  

The idea to define topological Hochschild homology by mimicking algebra in
this way is due to Goodwillie.  But because a symmetric monoidal category
of spectra was not available until recently, one could not simply implement
this idea.  B\"okstedt was the first one to define topological Hochschild
homology by modifying this idea to work with certain rings up to homotopy. 
This original definition of topological Hochschild homology 
concerns the case when $k=S$.
We restate the definition of the simplicial spectrum
$\thh_\point(R)$ and its realization, $\thh(R)$, from [B] 
for a symmetric ring spectrum.  See Definition \ref{thhr}.
In Theorem \ref{thm-thh-many} we show that for $k=S$ 
our new definitions are stably equivalent to the original definition
when $R$ is a cofibrant symmetric ring spectrum. 
As a corollary to this comparison theorem we see that B\"okstedt's 
definition of THH takes stable equivalences of $S$-algebras to
\p-isomorphisms.  Hence it always determines the right homotopy type, even
on non-connective and non-convergent ring spectra,
whereas the other two definitions give the right homotopy type only on
cofibrant symmetric ring spectra.  

\subsection{Two definitions of relative topological Hochschild homology}
The first definition corresponds to 
the derived tensor product notion of algebraic Hochschild homology.
The second definition mimics the Hochschild complex from algebraic
Hochschild homology.  As we see in Theorem \ref{thh comparison}, these 
notions are stably equivalent when $M$ is a cofibrant $R^e$-module.  

\begin{definition}\label{tk}Define $\tk(R;M)$ by $M \smasRe R$.  
\end{definition}

Let $\mu: R \smashk R  \xrightarrow{} R$ and $\eta: k \xrightarrow{} R$
be the multiplication and unit maps on $R$.  Let $\phi_r: M\smashk R 
\xrightarrow{} M$ and $\phi_l: R\smashk M \xrightarrow{} M$ be the right
and left $R$-module structure maps of $R$ acting on $M$. 
Let $R^{s}$ be the smash product over $k$ of $s$ copies of $R$, {\it i.e.,} 
$R\smashk \cdots
\smashk R$.  The following definition mimics the Hochschild complex
as in [CE]. 

\begin{definition}\label{tkt} $\tkt_{\point}(R;M)$ is the simplicial $k$-module with
$s$-simplices $M \smashk R^{s}$.  The simplicial face and degeneracy
maps are given by  
$$d_i= 
\begin{cases}
\phi_r \wedge (id_R)^{s-1} & \mbox {if } i=0 \\
(id_M) \wedge (id_R)^{i-1} \wedge \mu \wedge (id_R)^{s-i-1} & \mbox {if } 
1 \leq i < s\\
(\phi_l \wedge (id_R)^{s-1}) \circ \tau & \mbox{if } i=s
\end{cases}
$$
and $s_i= id_M \wedge (id_R)^i \wedge \eta \wedge (id_R)^{s-1}$. 
\end{definition}

Each level of this simplicial symmetric spectrum is a bisimplicial set.
Since the realization of bisimplicial sets is equivalent to taking
the diagonal, we use the diagonal to define the realization of this
simplicial symmetric spectrum.

\begin{definition}
Define the $k$-module $\tkt(R;M)$ as the diagonal of the bisimplicial set at
each level of this simplicial $k$-module.    
For the special cases of $k=S$ or $M=R$ we delete them from the notation. 
\end{definition}

Since the homotopy colimit of a diagram of symmetric spectra is
determined by the homotopy colimit of each level, the fact
that the homotopy colimit of a bisimplicial set is weakly
equivalent to the diagonal simplicial set, see [BK, XII 4.3], proves the
following proposition.

\begin{proposition}\label{diag=hocolim1}
The map $\hocolim_{\spectra}^{\Delta^{op}}\tkt_{\point}(R;M)
\xrightarrow{} \tkt(R;M)$ is a level equivalence.
\end{proposition}

Using $D$ and this Proposition we can show that  
the realization of a map which is a stable equivalence at each simplicial level
is a stable equivalence.

\begin{lemma}\label{stable hocolim}
Let $F, G \s B \to \spec$ be two diagrams of symmetric spectra with
a natural transformation $\eta \s F \to G$ between them.  If
$\eta(b) \s F(b) \to G(b)$ is a stable equivalence for each object 
$b$ in $B$ then $\hocolim^B F \to \hocolim^B G$ is a stable equivalence. 
\end{lemma}

\begin{proof}
Consider $D\eta \s DF \to DG$.  By Theorem \ref{detection-thm} this is a 
\p-isomorphism at each object, so by Proposition \ref{stable comparison}
the homotopy colimits are \p-isomorphic. 
Since $L'$, $F_0$ and homotopy colimits commute with homotopy colimits
and $\Omega^n$ commutes with homotopy colimits up to \p-isomorphism,
$\hocolim^B DF$ is \p-isomorphic to $D\hocolim^B F$.  Hence, $D\hocolim ^B F 
\to D\hocolim^B G$ is a \p-isomorphism.  Thus, by Theorem \ref{detection-thm},
$\hocolim^B F \to \hocolim^B G$ is a stable equivalence.
\end{proof}

\begin{corollary}\label{stable realization}
A map between simplicial symmetric spectra which is a stable equivalence
on each level induces a stable equivalence on the realizations. 
\end{corollary}  
 
\begin{proof}
This just combines Lemma \ref{stable hocolim} and [BK, XII, 4.3]. 
\end{proof}

\begin{proposition}\label{tkt2 comparison}
Let $R \xrightarrow{} R'$  be a 
stable equivalence between $k$-algebras which are cofibrant as $k$-modules,  
$M$ a $R^e$-module, $N$ a $(R')^e$-module, and 
$M \xrightarrow{} N$ a stable equivalence of $R^e$-modules.  
Then $\tkt(R;M) \xrightarrow{} \tkt(R';N)$ is a stable
equivalence.  In particular, $\tkt(R) \xrightarrow{} \tkt(R')$,
$\tkt(R;M) \to \tkt(R;N)$, 
and $\tkt(R;N) \to \tkt(R';N)$ are stable equivalences. 
\end{proposition}

First note that a cofibrant $k$-algebra is also cofibrant as a $k$-module
by Theorem \ref{R-alg}, so there are many examples of $k$-algebras which
are cofibrant as $k$-modules.  

\begin{proof}
Lemma \ref{R-cof.^stable=stable} applied
to $k$ shows that $P \smashk -$
preserves stable equivalences of $k$-modules if $P$ is a cofibrant
$k$-module.  Hence, $R^s \to R'{}^s$ is a stable equivalence between cofibrant
$k$-modules.  So both $M \smashk R^s \to N \smashk R^s$ and 
$N \smashk R^s \to N \smashk R'{}^s$ are also stable equivalences.  
Thus each simplicial level is a
stable equivalence.  Then Corollary \ref{stable realization} shows that this 
map induces a stable equivalence on $\tkt$.
\end{proof}

To compare these two definitions of topological Hochschild homology we 
first define certain
bar constructions.  Let $N$ be a left $R$-module, with $\phi_N\s R \smashk N \xrightarrow{} N$, and $M$ a right $R$-module, 
with $\phi_M\s M \smashk R \xrightarrow{} M$.
We define the topological bar construction 
$B^k_{\point}(M, R, N)$ by mimicking algebra. 

\begin{definition} The bar construction $B^k_{\point}(M,R,N)$
is the simplicial $k$-module with $s$-simplices $M \smashk R^{s} \smashk
N$.  The face and degeneracy maps are given by
$$d_i= 
\begin{cases}
\phi_M \wedge (id_R)^{s-1} \wedge id_N & \mbox{if } i=0\\
id_M \wedge (id_R)^{i-1} \wedge \mu \wedge (id_R)^{s-i-1} \wedge id_N 
& \mbox {if } 1 \leq i < s\\
id_M \wedge (id_R)^{s-1} \wedge \phi_N & \mbox{if } i=s
\end{cases}
$$
Let $B^k(M,R,N)$ be the realization of this simplicial $k$-module.
\end{definition} 

Let $c_{\point}(X)$ be the constant simplicial object with $X$ in each 
simplicial degree.  Using the identification $M \smasR R \iso M$, the map  
$\eta: k \xrightarrow{} R$ induces a simplicial $k$-module map 
$B^k_{\point}(M, R, N) \xrightarrow{} c_{\point}(M\smasR N)$.  

\begin{lemma}\label{bar} 
For $M$ a cofibrant $R$-module,  
the simplicial map of $k$-modules, $B^k_{\point}(M, R, N)
\xrightarrow{} c_{\point}(M\smasR N)$, induces a stable equivalence of
$B^k(M,R,N) \xrightarrow{} M\smasR N$. 
\end{lemma}

\begin{proof}
Note that $B^k_{\point}(M, R, N) \iso c_{\point}M \smasR B^k\point(R,R, N)$.
Since realization commutes with smash products, 
$B^k(M, R, N) \iso M \smasR B^k(R,R,N)$.  So using Lemma 
\ref{R-cof.^stable=stable} it is enough to show that $B^k(R,R,N) \xrightarrow{}
N$ is a stable equivalence.  The map $N \iso k \smashk N \xrightarrow{}
R\smashk N$ provides a simplicial retraction for $B^k_{\point}(R,R,N)$.
Hence the spectral sequence for computing the classical stable homotopy groups of
the homotopy colimit of this simplicial $k$-module collapses.  So the map
$B^k_{\point}(R,R,N) \xrightarrow{} c_{\point} N$ induces a 
\p-isomorphism on the realizations.  
\end{proof}   

Using the bar construction we now show that the two definitions of
topological Hochschild homology are stably equivalent when $M$ is
a cofibrant $R^e$-module.

\begin{theorem} \label{thh comparison}
There is a natural map of $k$-modules 
$\tkt(R;M) \xrightarrow{} \tk(R;M)$ which is a 
stable equivalence for $M$ a cofibrant $R^e$-module. 
\end{theorem}

\begin{proof}
We show that $\tkt(R;M)$ is naturally isomorphic to 
$M \smasRe B^k(R,R,R)$ below.  Then the map $\tkt(R;M) \xrightarrow{}
\tk(R;M)$ is given by $M\smasRe \phi$ for 
$\phi:B^k(R,R,R) \xrightarrow{} R$.  
$R$ is always a cofibrant $R$-module,
hence $\phi$ is a 
stable equivalence by Lemma \ref{bar}.  Then Proposition 
\ref{R-cof.^stable=stable} shows that
$M \smasRe \phi$ is a stable equivalence since $M$ is a cofibrant
$R^e$-module.

To see that $\tkt(R;M)$ is naturally isomorphic to $M \smasRe B^k(R,R,R)$
we show that $\tkt_{\point}(R;M)$ is naturally isomorphic to $c.(M) 
\smasRe B^k_{\point}(R,R,R)$.  On each simplicial level
there are natural isomorphisms $$ M\smashk R^{s} \iso M \smasRe
(R^e \smashk R^{s}) \iso M \smasRe (R \smashk R^{s} \smashk R)
= M \smasRe B^k_s(R,R,R).$$  These isomorphisms commute with the
simplicial structure.  Hence the simplicial $k$-modules are naturally 
isomorphic, so their realizations are also naturally isomorphic. 
\end{proof}

\subsection{B\"okstedt's definition of topological Hochschild homology}
 
We now define the simplicial spectrum $\thh_\point (R;M)$ and
its realization $\thh(R;M)$ following B\"okstedt's original definitions.
Each of the levels of the simplicial spectrum $\thh_\point$ can
be defined for a general symmetric spectrum $X$.  A ring structure is
only necessary for defining the simplicial structure.  In fact,
each of the levels of $\thh_\point$ can be thought of as a functor which 
gives the correct \p-isomorphism type for the smash product of 
symmetric spectra.   We start by considering each of these levels as
a functor of several variables.  

Let ${\bf X}$ denote a sequence of $j+1$ spectra, $X^0, \dots, X^j$. 
Define a functor $\D^j{\bf X}$ from $I^{j+1}$ to $\spec$  which
at ${\bf n}= (n_0, \cdots, n_j)$ takes the value,
$${\mathcal D}^j{\bf X}({\bf n})=\Omega^n L'F_0(X^0_{n_0} \sm \dots \sm
 X^j_{n_j})$$
where $ L'$ is a level fibrant replacement functor and 
$n=\Sigma n_i$, the sum of the $n_i$.
Note that ${\mathcal D}^0(X) $ is ${\mathcal D}_X$, the functor defined
at the beginning of section \ref{D}.  To see   that  ${\mathcal D}^j
{\bf X}$ is defined over $I^{j+1}$ one  uses maps similar to those
described for ${\mathcal D}_X$.

\begin{definition}\label{defn-thh_j-many}
Let $X^0, \dots, X^j$  be symmetric spectra.     Define
$$\h_j{\bf X}= \hocolim^{I^{j+1}}{\mathcal D}^j{\bf X}.$$
\end{definition}

 We now define a
natural transformation $\phi_j{\bf X}\s \h_j{\bf X}
\to D(X^0 \sms \dots \sms X^j)$.
Let $\mu:I^{j+1} \xrightarrow{} I$ be the functor induced by concatenation
of all of the factors in $I^{j+1}$.  Then there is a natural transformation
from ${\mathcal D}^j{\bf X}$ to
$\mu^*{\mathcal D}^0(X^0 \sms \dots \sms X^j)$.  This natural
transformation is induced by the map from $X^0_{n_0} \wedge \cdots \sm
X^j_{n_j}$ to  the $n$th level of $X^0 \sms \dots \sms X^j.$
This map is
$\Sigma_{n_0} \times \cdots \times \Sigma{n_j}$ equivariant, which
is exactly what is necessary over $I^{j+1}$.
Hence, on homotopy colimits there is a natural map
$\hocolim^{I^{j+1}}{\mathcal D}^j{\bf X}
\xrightarrow{} \hocolim^{I^{j+1}}\mu^*{\mathcal D}^0(X^0 \sms \dots \sms X^j)$.

\begin{definition}
There is a natural transformation
$\phi_j{\bf X}\s \h_j{\bf X}
\to D(X^0 \sms \dots \sms X^j)$.  It is given by the composition
\[\hocolim^{I^{j+1}}{\mathcal D}^j{\bf X}
\xrightarrow{} \hocolim^{I^{j+1}}\mu^*{\mathcal D}^0(X^0 \sms \dots \sms X^j)
\]
\[
\xrightarrow{} \hocolim^{I}{\mathcal D}^0(X^0 \sms \dots \sms X^j).\]
\end{definition}

\begin{proposition}\label{prop-phi_j-many}
For any cofibrant symmetric spectra, $X^0, \dots, X^j$, the map
$\phi_j{\bf X}$ is a \p-isomorphism.
\end{proposition}

This proposition is proved in subsection \ref{sub-prop}.  It is used in 
proving the comparison theorem between B\"okstedt's definition of $\THH$ and
our previous definition of $\tHH$.  As a corollary of this proposition, 
$\h_j$ gives the correct \p-isomorphism type for the derived smash product of 
$j+1$ symmetric spectra.  Recall that the smash product is only homotopy
invariant on cofibrant spectra, so the derived smash product is the smash
product of the cofibrant replacements. 
In the stable model category of symmetric spectra, consider a cofibrant 
replacement functor, $C$, analogous to the fibrant replacement functor $L$. 

\begin{corollary}
\p$\h_j{\bf X}$
is isomorphic to \p$L(CX^0 \sms \dots \sms CX^j)$, the derived homotopy of
the derived smash product of $X^0,  \dots, X^j$.  
\end{corollary}

\begin{proof}
Since $C$ is a cofibrant replacement functor, $CX \to X$ is a level 
equivalence.  Hence
$\h_j(CX^0, \dots, CX^j) \to \h_j(X^0, \dots, X^j)$ is a level
equivalence by Proposition \ref{level comparison} because the map is a level 
equivalence at each object in the diagram defining $\h_j$.  So this
corollary follows from  Proposition \ref{prop-phi_j-many} since \p$D(CX^0 \sms
\dots \sms CX^j)$ is isomorphic to \p$L(CX^0 \sms \dots \sms CX^j)$ by 
Theorem \ref{thm-derived}. 
\end{proof}

We now use this functor $\h_j$ to define $\THH$ following B\"okstedt's
definition in [B]. 

\begin{definition}Let $R$ be a symmetric ring spectrum with $M$ an $R^e$-module.
Define $\THH_j(R;M)=\h_j(M, R, \dots, R)$.  
\end{definition}

The functors $\thh_j(R;M)$ fit together
to form a simplicial symmetric spectrum $\thh_{\textstyle \cdot}(R;M)$.
Although the definition of $\thh_j(R;M)$
does not use the ring structure of $R$ or the module structure of $M$, the 
simplicial structure of
$\thh_{\textstyle \cdot}(R;M)$ does use both the multiplication and unit maps.
The $i$th face map uses the
functor $\delta_i\s I^{j+1}\xrightarrow{} I^j$ defined by concatenation of
the sets in factors $i$ and $i+1$.  The last face map uses the 
cyclic permutation of
$I^{j+1}$ followed by concatenation of the first two factors.
For ease of notation let ${\mathcal D}^j(R;M)= {\mathcal D}^j(M, R, \dots, R)$.
The multiplication of $R$ and $M$ defines a natural transformation of functors
from ${\mathcal D}^j(R;M)$ to 
$\delta_i^*{\mathcal D}^{j-1}(R;M)$.
So $d_i$ is the composition 
$$d_i: \hocolim^{I^{j+1}}{\mathcal D}^{j}(R;M) \xrightarrow{} 
\hocolim^{I^{j+1}}\delta_i^*{\mathcal D}^{j-1}(R;M)
\xrightarrow{} \hocolim^{I^j}{\mathcal D}^{j-1}(R;M).$$
The degeneracy maps are similar.

\begin{definition}\label{thhr}
Define $\THH(R;M)$ as the diagonal of the bisimplicial set at each
level of the simplicial symmetric spectrum $\THH_{\point}(R;M)$.
\end{definition}

One can check that each level in this spectrum agrees with the
definition in [B] when $M=R$. 

As in Proposition \ref{diag=hocolim1} we have the following equivalence. 

\begin{proposition}\label{diag=hocolim}
The map $\hocolim_{\spectra}^{\Delta^{op}}\thh_{\point}(R;M) 
\xrightarrow{} \thh(R;M)$ is a level equivalence.
\end{proposition}

The next theorem shows that the definition of topological Hochschild
homology which mimics the Hochschild complex is stably equivalent
to the original definition of topological Hochschild homology.  

\begin{theorem}\label{thm-thh-many}
Let $R$ be a cofibrant ring spectrum.  Then there is a natural
zig-zag of stable equivalences between $\tHH(R;M)$ and $\THH(R;M)$.
\end{theorem}

\begin{proof}
The zig-zag of functors between $\tHH$ and $\THH$ is induced by a  
zig-zag of maps between the simplicial complexes defining $\tHH$ and $\THH$.
First one applies the zig-zag of functors 
$1 \xrightarrow{\psi} L \xrightarrow{} ML  \xleftarrow{} DL \xleftarrow{D 
\psi} D$ to each simplicial level of the Hochschild complex defining $\tHH$.  
Here, $L$ is the fibrant
replacement functor, $M$, $D$, and the natural transformations are defined in 
section \ref{D}, see \ref{def D}, \ref{def M}, and the proof of 
\ref{detection-prop} part 3.   Then there is a
natural map $\phi_j \s \THH{}_j(R;M) \xrightarrow{} D(M \sms R^{j})$ 
To see that the $\phi_j$ maps commute with the simplicial maps, one needs to
note that the multiplication maps commute with the first map in the composite
defining $\phi_j$.  This follows since the map $R_n \wedge R_m \to R_{n+m}$
is the map on the appropriate wedge summand of the map 
$R\otimes R \to R$ which induces the map $R \wedge_S R \to R$.  The maps 
involving $M$ are similar. 
Putting these simplicial levels together one gets
a zig-zag of natural transformations from $\tHH{}_\point(-)$ to 
$\THH_\point(-)$.  Note that by composing maps this zig-zag is
only of length 2. 

The zig-zag of functors between $1$ and $D$ was investigated in section 
\ref{D}.  $\psi$ induces a stable equivalence on any object by definition
of the fibrant replacement functor $L$.  
So Corollary \ref{stable realization}
applies to show that $\psi$ induces a stable equivalence on the
realization of these simplicial $S$-modules.  

Since $\psi$ is always a stable equivalence, $D \psi$ is a
\p-isomorphism on each simplicial level by Theorem 
\ref{detection-thm}.  The two natural transformations of 
middle functors $L \xrightarrow{} ML \xleftarrow{} DL$ induce level 
equivalences by Proposition \ref{detection-prop}.  
Propositions \ref{level
comparison} and \ref{stable comparison} apply to these natural 
transformations to show that they also induce level equivalences and
\p-isomorphisms on the realizations. 

So the only part of the zig-zag between $\tHH(R;M)$ and $\THH(R;M)$ that is left
is $\THH_{\point}(R;M) \xrightarrow{} D(\tHH_{\point}(R;M))$.    
Let $CM \to M$ be a cofibrant replacement of $M$ as an $R^e$-module.  
Then by Proposition~\ref{tkt2 comparison}, $\tHH_j(R;CM) \to \tHH_j(R;M)$
is a stable equivalence.  Similarly, $\THH_j(R;CM)\to \THH_j(R;M)$ is a
stable equivalence since $CM \to M$ is a level equivalence and hence induces
a level equivalence on the homotopy colimits used to define $\THH_j$.
So it is enough to prove that $\THH_{\point}(R;M) \xrightarrow{} D(\tHH_{\point}(R;M))$ is a stable equivalence in the case when $M$ is cofibrant.    

Since $R$ is cofibrant as an $S$-algebra, it is also cofibrant as an $S$-module.
Since $M$ is cofibrant as an $R^e$-module and $R^e$ is cofibrant, $M$ is also
cofibrant as an $S$-module.
Proposition \ref{prop-phi_j-many} shows that if $R$ and $M$ are any cofibrant
$S$-modules then $\THH{}_j(R;M) \xrightarrow{} D(M \sms R^{j})$ is a 
\p-isomorphism.  Then Proposition \ref{stable comparison} shows that this
is enough to ensure that the map on the realizations 
is a \p-isomorphism.  Hence, assuming Proposition \ref{prop-phi_j-many}, 
this finishes the proof of Theorem \ref{thm-thh-many}.  
\end{proof}

\begin{corollary}
The derived stable homotopy groups of $\tHH(R;M)$ are isomorphic to
\p$\THH(R;M)$. 
\end{corollary}

Since $\tHH(R;M)$ and $\THH(R;M)$ are stably homotopic their derived stable
homotopy groups must be isomorphic.  So this corollary says
that the derived stable homotopy groups of $\tHH(R;M)$ are isomorphic to
the classical stable homotopy groups of $\THH(R;M)$.

\begin{proof}
The proof of Theorem \ref{thm-thh-many} shows that the map from $\THH(R;M)$
to the realization of $D\tHH_{\point}(R;M)$ is a \p-isomorphism.  
But this realization is \p-isomorphic to $D\tHH(R;M)$ 
as shown in the proof of Proposition \ref{stable hocolim}.  
So the derived stable homotopy groups of $\tHH(R;M)$, {\it i.e.,} \p$D\tHH(R;M)$,
are isomorphic to the classical stable homotopy groups of $\THH(R;M)$.
\end{proof}

Using this comparison we can show that B\"okstedt's original definition
of $\THH$ takes stable equivalences of ring spectra to
\p-isomorphisms.  This is a stronger result than for $\tHH$ because
no cofibrancy condition is needed here and the map is a \p-isomorphism,
not just a stable equivalence.   

\begin{corollary}\label{thh-stable}
Let $R \xrightarrow{} R'$ be a stable equivalence of ring spectra, $M$
a $R^e$-module, $N$ a $(R')^e$-module, and $M \to N$ a stable equivalence
of $R^e$-modules.  
Then $\thh(R;M) \xrightarrow{} \thh(R';N)$ is a \p-isomorphism.
\end{corollary}

\begin{remark}
This corollary could also be proved without using these comparison  
results.  One can show that each $\thh_j$ takes stable equivalences
to \p-isomorphisms following arguments similar to those for
$\thh_0=D$ in section \ref{D}.  Then Proposition \ref{stable comparison}
shows that the realization, $\thh$, also takes stable equivalences
to \p-isomorphisms.  
\end{remark}
 
\begin{proof}
In the category of symmetric ring spectra, define a functorial 
cofibrant replacement functor, $C$.   
Applying this functor we have the following square.
\[
\begin{CD}
CR @>>> CR'\\
@VVV   @VVV\\
R @>>> R'
\end{CD}
\]
Each of the vertical maps is a level trivial fibration and hence a level 
equivalence.  The bottom map is a stable equivalence by assumption.  Hence
the top map is also a stable equivalence.  To show that $\thh$ applied
to the bottom map is a \p-isomorphism we show that $\thh$ applied
to the other three maps in this square are \p-isomorphisms.  

We also need to consider cofibrant replacements of the modules in question.
$M$ is a $(CR)^e$-module and $N$ is a $(CR')^e$-module.  We replace
them by modules which are cofibrant as underlying $S$-modules.  Since
$CR$ is a cofibrant $S$-algebra it is a cofibrant $S$-module.  Thus $(CR)^e$ 
is also cofibrant as an $S$-module by Proposition \ref{monoidal}.  Hence
the cofibrations in the category of $(CR)^e$-modules are also underlying
cofibrations.  So let $CM \to M$ be the cofibrant replacement of $M$ in
the category of $(CR)^e$-modules.  Similarly, $R'$ 
let $CN \to N$ be the cofibrant replacement of $N$ as a $(CR')^e$-module.  
Then both $CM$ and $CN$ are cofibrant as $S$-modules.
Also, in the category of $(CR)^e$-modules by the lifting property in the model 
category of $(CR)^e$-modules we have a map $CM \to CN$ because
$CN \to N$ is a level trivial fibration.  This map $CM \to CN$ is a stable
equivalence by the two out of three property.  

The level equivalences $CR \to R$ and $CM \to M$ induce a level equivalence on 
each object of the diagram defining $\thh_j$.  So by applying Proposition 
\ref{level comparison} and Lemma \ref{diag=hocolim} this shows that 
$\thh(CR; CM) \to \thh(R;M)$ is a level equivalence.  Similarly
$\thh(CR'; CN) \to \thh(R';N)$ is a level equivalence.  

For the top map, first consider applying $\tHH$.  Proposition 
\ref{tkt2 comparison} implies that $\tHH(CR;CM) \xrightarrow{} \tHH(CR';CN)$ 
is a stable equivalence.
Hence by Theorem \ref{detection-thm}, $D\tHH(CR;CM) \xrightarrow{} 
D\tHH(CR';CN)$ is a \p-isomorphism.  But in the proof of Theorem 
\ref{thm-thh-many}  we showed that $\THH \xrightarrow{} D\tHH$
induces a \p-isomorphism if the ring and module are cofibrant as $S$-modules.  
So $\THH(CR;CM) \xrightarrow{} \THH(CR';CN)$ is a \p-isomorphism.  
Stringing these equivalences together finishes the proof of
this corollary.  Because Proposition \ref{prop-phi_j-many} applies to each
level, we have actually shown that each 
$\THH_j(R;M) \xrightarrow{} \THH_j(R';N)$ is also a \p-isomorphism. 
\end{proof}

\subsection{Proof of Proposition \ref{prop-phi_j-many}}\label{sub-prop}
To prove Proposition \ref{prop-phi_j-many} we follow an outline similar to the
proof that $D$ takes stable trivial cofibrations to \p-isomorphisms.  We show
that $\phi_j$ is a \p-isomorphism when it is evaluated only on
free symmetric spectra, {\it i.e.,} some $F_nK$.  
Then we prove an induction step lemma which
deals with pushouts over generating stable trivial cofibrations.  Using
these lemmas we show that $\phi_j$ is a \p-isomorphism on
any collection of cofibrant spectra.   

\begin{lemma}\label{lem-free}
$\phi_j(F_{n_0}K_0, \dots, F_{n_j}K_j)$ is a \p-isomorphism.
\end{lemma}

\begin{lemma}\label{lem-po-many}
Let $A \to B$ be a stable cofibration and $X^0, \dots, X^j$ be cofibrant
$S$-modules.  Consider the following pushout square.
\[
\begin{CD}
A @>>> X \\
@VVV @VVV\\ 
B @>>> Y
\end{CD}
\]
Assume that $\THH_{j+1}(X^0, \dots, Z, \dots, X^j) \to D(X^0 \sms \dots \sms
Z \sms \dots \sms X^j)$ is a \p-isomorphism for $Z = A, B,$ or $X$ where
$Z$ is inserted between the $i$th and $i+1$st spots.   Then
$\THH_{j+1}(X^0, \dots, Y,\dots, X^j) \to D(X^0 \sms \dots \sms Y \sms \dots
\sms X^j)$ is a \p-isomorphism.
\end{lemma}

Using these two lemmas we can now prove Proposition \ref{prop-phi_j-many}.

\begin{proof}[Proof of Proposition \ref{prop-phi_j-many}]
We prove this by induction on $i$ with the induction assumption that $\phi_j$
is a \p-isomorphism when $j-i$ variables are free spectra and the
other variables are cofibrant.  Lemma \ref{lem-free} verifies this for
$i=0$.  For the induction step, in one variable 
we build up a cofibrant spectrum from the initial
spectrum by retracts, colimits, and pushouts over generating cofibrations.
Since retracts of \p-isomorphisms are \p-isomorphisms and
$\phi_j$ of a retract is a retract we only need to consider colimits
and pushouts.
 
Because $F_0$, smash products, $L'$, $\Omega^n$, and homotopy colimits commute 
with filtered colimits, $\h_j$ of a colimit in one of the variables 
is a colimit.  This is also true of $D$.  Since a filtered colimit of
\p-isomorphisms is a \p-isomorphism this means $\phi_j$ of a colimit in
one variable is a \p-isomorphism if it is a \p-isomorphism at each spot
in the sequence.  Hence we are only left with pushouts.  

Since $\phi_j$ is a level equivalence between trivial spectra if one of
the variables is the initial spectrum, $*$, one can proceed by induction
to verify the pushout property.  By induction the two corners in the
pushout corresponding to the generating cofibration are \p-isomorphisms. This 
is because generating cofibrations are of the form $F_nK \to F_nL$, so these 
two corners have one extra variable a free spectrum
and hence fall into the case covered by the previous induction step.
The third corner is assumed to be a \p-isomorphism
by induction, hence $\phi_j$ is a \p-isomorphism on the pushout corner by
Lemma \ref{lem-po-many}.   
\end{proof}

\begin{proof}[Proof of Lemma \ref{lem-free}]
To show that $\phi_j(F_{n_0} K_0 , \dots , F_{n_j}K_j)$ is a \p-isomorphism
we first establish the stable homotopy type of $\h_j(F_{n_0} K_0 , \dots , 
F_{n_j}K_j)$.  There is a functor 
$\F_{(n_0, \dots, n_j)}X \s I^{j+1} \to \spec$
defined by $$\F_{(n_0, \dots, n_j)}X (m_0, \dots m_j)=
\hom_{I^{j+1}}((n_0, \dots, n_j),(m_0, \dots m_j))_+ \sm X.$$ Then $\F_{(n_0, 
\dots, n_j)}(-)$ is left adjoint to the functor from I-diagrams over $\spec$ to
$\spec$ which evaluates the diagram at $(n_0, \dots, n_j)\in I^{j+1}$.
There is a map of diagrams $\F_{(n_0, \dots, n_j)}(\Omega^{ n} 
L'F_0(K_0 \sm \dots \sm K_j)) \to \D^j(F_{n_0} K_0 , \dots , F_{n_j}K_j)$
where $n=\Sigma n_i$.     
Each spot in this diagram is a \p-isomorphism.  This is similar to
the proof of Lemma \ref{stable type}, on each level the map is an equivalence
in the stable range by the Blakers-Massey and the Freudenthal
suspension theorems.  Hence the map on homotopy colimits is also a 
\p-isomorphism, $\Omega^{ n}L' F_0(K_0 \sm \dots \sm K_j)) \to \h_j
(F_{n_0} K_0 , \dots , F_{n_j}K_j)$.    

By Lemma \ref{stable type}, $\Omega^{n} L'F_0(K_0 \sm \dots \sm K_j)) \to
D(F_{n} (K_0 \sm \dots 
\sm K_j))$ is also a \p-isomorphism.  To see that $\phi_j$ induces a 
\p-isomorphism, note that on the free diagrams there are similar maps
$\hocolim^{I^{j+1}}\F_{(n_0, \dots, n_j)}
(\Omega^{n}L' F_0(K_0 \sm \dots \sm K_j))  \to \hocolim^{I^{j+1}}
\mu^*\F_{ n} (\Omega^{n}L' F_0(K_0 \sm \dots \sm K_j)) \to  
\hocolim^{I} \F_{ n} (\Omega^{ n}L' F_0(K_0 \sm \dots \sm K_j))$
which induce level equivalences on the homotopy colimits.  
\end{proof}

To prove Lemma \ref{lem-po-many} we first need to show that $\h_j$ of
a homotopy pushout in one variable is a homotopy pushout.  

\begin{lemma}\label{THH-p.o.}
Let $X^0, \dots, X^{j}$ be cofibrant spectra.  If
\[
\begin{CD}
A @>>> X\\
@VVV   @VVV\\
B @>>> Y
\end {CD}
\]
is a pushout square with $A \xrightarrow{} B$ a cofibration, then
\[
\begin{CD}
\h_{j+1}(X^0, \dots, A, \dots X^{j})@>>>\h_{j+1}(X^0, \dots, X, \dots X^{j}) \\
@VVV   @VVV\\
\h_{j+1}(X^0, \dots,B, \dots X^{j})@>>>\h_{j+1}(X^0, \dots, Y, \dots X^{j}) 
\end {CD}
\]
is a homotopy pushout square.  I.e, if $P$ is the homotopy pushout 
of the second square then $P \xrightarrow{} 
\h_{j+1}(X^0, \dots, Y, \dots X^{j}) $ is a stable equivalence.
In fact, $P \xrightarrow{}\h_{j+1}(X^0,\dots, Y,\dots X^{j})$ is a 
\p-isomorphism.
\end{lemma}

\begin{proof}
This proof is similar to the proof of Lemma \ref{D-p.o.}.
As with Lemma \ref{D-p.o.}, it is enough to consider each object in
$I^{j+1}$ since homotopy colimits commute. 

The following square is a pushout square with the left map a cofibration
since each $X^i$ is cofibrant.
\[
\begin{CD}
X^0_{n_0} \sm \dots \sm A_{n_i} \sm \dots X^{j}_{n_j} @>>> 
X^0_{n_0} \sm \dots \sm X_{n_i} \sm \dots X^{j}_{n_j}\\ 
 @VVV @VVV\\ 
X^0_{n_0} \sm \dots \sm B_{n_i} \sm \dots X^{j}_{n_j} @>>> 
X^0_{n_0} \sm \dots \sm Y_{n_i} \sm \dots X^{j}_{n_j} 
\end{CD}
\]

The first step in constructing $\h_j$ is just applying $F_0$ to this square.
$F_0$ preserves cofibrations and pushouts, hence $F_0$ applied to this
square is a homotopy pushout.  
$L'$ preserves homotopy pushout squares up to level equivalence and
$\Omega^{\Sigma n_i}$ preserves homotopy pushout squares up to \p-isomorphism.
Hence the map from the homotopy pushout
to the bottom right corner is a \p-isomorphism.  Since the homotopy colimit of
\p-isomorphisms is a \p-isomorphism, this finishes the proof.
\end{proof}

\begin{proof}[Proof of Lemma \ref{lem-po-many}]
Both $\h_j$ and $D$ take homotopy pushouts in one variable to
homotopy pushouts where the map from the pushout to the bottom right corner
is a \p-isomorphism by Lemmas \ref{THH-p.o.} and \ref{D-p.o.}.  Hence, this 
lemma follows from the fact that homotopy colimits
preserve \p-isomorphisms,  Lemma \ref{stable comparison}.    
\end{proof}

\nocite{Bo}
\nocite{BF}
\nocite{BK}
\nocite{MMSS}
\nocite{HM}
\nocite{HSS}
\nocite{SS}
\nocite{F1}
\nocite{F2}
\nocite{DS}
\nocite{S}
\nocite{CE}
\nocite{Q}
\bibliography{a}
\bibliographystyle{plain}

\bibliographystyle{amsalpha}
\end{document}